\documentclass[10pt,twocolumn,twoside]{IEEEtran}

\usepackage{savesym}
\savesymbol{proof}
\savesymbol{endproof}
\usepackage{amsmath,amsthm,amssymb,amsfonts}
\restoresymbol{AMS}{proof}
\restoresymbol{AMS}{endproof}
\usepackage{cite}
\usepackage{graphicx}
\graphicspath{ {img/} }
\usepackage{algorithm}
\usepackage{algpseudocode}

\newtheorem{definition}{Definition}
\newtheorem{theorem}{Theorem}
\newtheorem{proposition}{Proposition}
\newtheorem{lemma}{Lemma}
\newtheorem{corollary}{Corollary}
\newtheorem{remark}{Remark}
\newtheorem{assumption}{Assumption}
\newtheorem{problem}{Problem}

\usepackage{enumerate}

\usepackage[svgnames]{xcolor}
\usepackage{tikz}
\usetikzlibrary{arrows}
\usepackage{verbatim}
\usetikzlibrary{arrows,shapes,automata}
\usetikzlibrary{fit}

\tikzset{
	treenode/.style = {align=center, inner sep=0pt, text centered,
		font=\sffamily},
	task/.style = {treenode, circle, white, font=\sffamily\bfseries, draw=black,
		fill=black, text width=1.0em},
	agent/.style = {treenode, circle, SteelBlue, draw=SteelBlue,
		text width=1.0em, thick},
}
\usepackage{textcomp}
\begin{document}
\title{A Distributed Augmenting Path Approach for the Bottleneck Assignment Problem}
\author{Mitchell Khoo, Tony A. Wood, Chris Manzie and Iman Shames
\thanks{The research is funded by Defence Science and Technology Group through research agreements MyIP: 7558 and MyIP: 7562. }
\thanks{All authors are with the Department of Electrical and Electronic Engineering at the University of Melbourne, (e-mails: khoom1@student.unimeld.edu.au, \{wood.t,manziec,iman.shames\}@unimelb.edu.au). }}

\maketitle

\begin{abstract}
We develop an algorithm to solve the Bottleneck Assignment Problem (BAP) that is amenable to having computation distributed over a network of agents. This consists of exploring how each component of the algorithm can be distributed, with a focus on one component in particular, i.e., the function to search for an augmenting path. An augmenting path is a common tool used in most BAP algorithms and poses a particular challenge for this distributed approach. Given this significance, we compare two different methods to search for an augmenting path in a bipartite graph. We also exploit properties of the augmenting paths to formalise conditions for which the solution from subsets of the sets of agents and tasks can be used to solve the BAP with the full sets of agents and tasks. In the end, we evaluate and compare the derived approaches with a numerical analysis.
\end{abstract}


\section{Introduction} \label{sec:introduction}
An assignment problem emerges when multiple tasks must be allocated to multiple agents. In \cite{taxonomy,assignprob,survey}, reviews of the different types of assignment problems are presented and categorised according to objective functions, the number of agents required to complete tasks and the number of tasks individual agents can carry out. The bottleneck assignment problem (BAP) is a particular assignment problem with the objective of minimising the costliest agent-task pairing. Such an assignment is typically necessary in time-critical problems, where tasks are carried out by agents simultaneously and the time taken to complete all tasks must be minimised. An example application of the BAP is the threat seduction problem in \cite{decoys}, where decoys are assigned to multiple incoming threats.

The BAP is well-studied and there exist many centralised algorithms to solve it, e.g., \cite{decoys,cBAP2,cBAP3,cBAP4,cBAP5}. Centralised computation refers to the case where one decision-maker has access to the allocation cost of every potential agent-task pairing. The threshold algorithm in \cite{cBAP2} involves an initial threshold cost and checking if an assignment can be created with allocation costs smaller than the threshold. The threshold is iteratively increased until a valid assignment is found, which then corresponds to the optimal assignment for the BAP. In \cite{cBAP3,cBAP4}, improvements on the complexity of the threshold algorithm are made by moving the threshold according to a binary search pattern rather than incrementally. 
In~\cite{cBAP5}, the BAP is solved for a subset agents and tasks, where the size of the subset is increased until it contains all the agents and tasks. 



In many applications, no centralised decision-maker is available and there are restrictions on the information shared between agents. This is motivation to consider solving the BAP with computation distributed over a network of agents. As a comparison, there are many distributed algorithms for assignment problems with other objectives. For example, in the linear assignment problem (LAP), tasks are allocated to agents such that the sum of the costs of the allocations is minimised. 
In~\cite{distHun}, a distributed version of the well-known Hungarian Method from~\cite{hungarian} for solving the LAP is developed. In~\cite{auc1,auc2}, so-called auction algorithms are presented. The LAP can be cast as a linear program with binary decision variables to represent the allocation of an agent to a task. In~\cite{simplex}, a distributed simplex algorithm is developed to solve linear programs and the LAP is used as a motivating example.
A greedy algorithm sequentially picks one allocation of a task to an agent with lowest cost from the remaining choices of allocations. There exist distributed greedy algorithms, e.g., the Consensus-Based Auction Algorithm (CBAA) in \cite{CBAA}.  

We present a distributed algorithm for solving the BAP based on the framework introduced in~\cite{distBAP}. We apply tools introduced in~\cite{structure} to analyse the graph theoretical structure of the problem and propose different methods for implementing the algorithm efficiently. The concept of verifying the existance of an augmenting path plays a crucial role in solving the BAP. Thus, one of the main contributions of this paper is to explore two different distributed methods to search for an augmenting path in a bipartite graph and to compare the trade-offs between them. We also consider the situation where the set of agents and tasks are partitioned, i.e., the BAP is split into two smaller problems. We derive conditions for which a solution of the overall problem can be found by exploiting the solutions of the two smaller problems.

The rest of this paper is organised as follows. In Section~\ref{sec:prelim}, the graph theoretical background is introduced, including several novel definitions for analysing the structure of the BAP. In Section~\ref{sec:dist}, we present an algorithm for solving the BAP with computation distributed over a network of agents, which includes a distributed search for an augmenting path. In Section~\ref{sec:alternative}, we present an alternative approach for conducting the distributed search for an augmenting path and compare it to the approach in Section~\ref{sec:dist}. In Section~\ref{sec:exploit}, we explore how the partitioning of the network of agents can be exploited to solve the BAP. In Section~\ref{sec:numerical}, we present a numerical analysis of the proposed approaches.


\section{Preliminaries} \label{sec:prelim}

Consider an undirected graph $\mathcal{G}=(\mathcal{V,E})$, where $\mathcal{V}$ is a set of vertices and $\mathcal{E}$ is a set of edges. We provide the following graph theoretical definitions, which are found in~\cite{assignprob,graph}.

\begin{definition}[Matching, Maximum Cardinality Matching] \label{def:matching}
	A matching $\mathcal{M}$ of graph $\mathcal{G}$ is a set of edges such that $\mathcal{M}\subseteq \mathcal{E}$ and no vertex $v\in\mathcal{V}$ is incident to more than one edge in $\mathcal{M}$. A Maximum Cardinality Matching (MCM) is a matching $\mathcal{M}_{max}$ of $\mathcal{G}$ with maximum cardinality.
\end{definition}

\begin{definition}[Neighbours] \label{def:neigh}
	The set of neighbours of vertex $v\in\mathcal{V}$ in graph $\mathcal{G}$ is $N(\mathcal{G},v):=\{k|\{v,k\}\in \mathcal{E}\}$.
\end{definition}


\begin{definition}[Path] \label{def:path}
	Given distinct vertices $v_1, v_2, ..., v_{l+1}\in \mathcal{V}$ such that for all $k\in\{1,2,...,l\}$, $v_{k+1}\in N(\mathcal{G},v_k)$, the set of edges $P(\mathcal{E},v_1,v_{l+1}):=\{ \{v_k,v_{k+1}\}|k\in\{1,2,...,l\}\}$ is a path between $v_1$ and $v_{l+1}$ of length $l$.
\end{definition}

\begin{definition}
	(Diameter) The diameter of graph $\mathcal{G}$ is $D(\mathcal{G}):=\max_{v,v'\in \mathcal{V}} l_{vv'}$, where $l_{vv'}$ is the length of the shortest path between vertices $v,v'\in \mathcal{V}$.
\end{definition}

\begin{definition}[Alternating path] \label{def:alt}
	Given a matching $\mathcal{M}$ and a path $\mathcal{P}$, $\mathcal{P}$ is an alternating path relative to $\mathcal{M}$ if and only if each vertex that is incident to an edge in $\mathcal{P}$ is incident to no more than one edge in $\mathcal{P}\cap\mathcal{M}$ and no more than one edge in $\mathcal{P}\backslash\mathcal{M}$.
\end{definition}

If the elements of an alternating path relative to $\mathcal{M}$ are arranged in a sequence $\{v_1,v_2\}$, $\{v_2,v_3\}$, ..., $\{v_{l-1},v_l\}$, $\{v_l,v_{l+1}\}$, then the edges in the sequence alternate between edges in $\mathcal{M}$ and edges not in $\mathcal{M}$.

\begin{definition}[Free vertex] \label{def:free}
	Given a matching $\mathcal{M}$, vertex $v\in \mathcal{V}$ is free if and only if $\{v,w\}\notin \mathcal{M}$, for all $w\in \mathcal{V}$.
\end{definition}

\begin{definition}[Augmenting path] \label{def:aug}
	Given a matching $\mathcal{M}$ and a path $\mathcal{P}$ between vertices $v_1,v_{l+1}\in\mathcal{V}$, $\mathcal{P}$ is an augmenting path relative to $\mathcal{M}$ if and only if $\mathcal{P}$ is an alternating path relative to $\mathcal{M}$ and $v_1$ and $v_{l+1}$ are both free vertices.
\end{definition}

A tree is a graph $\mathcal{G}_t=(\mathcal{V}_t,\mathcal{E}_t)$ for which any two vertices in $\mathcal{V}_t$ are connected by exactly one path, and one vertex $r\in\mathcal{V}_t$ is designated as the root vertex.

\begin{definition}[Level] \label{def:level}
	Given a tree $\mathcal{G}_t=(\mathcal{V}_t,\mathcal{E}_t)$ with root vertex $r\in\mathcal{V}_t$, the level of a vertex $i\in\mathcal{V}_t$ in $\mathcal{G}_t$ is the length of the path between $i$ and $r$.
\end{definition}

\begin{definition}[Parent] \label{def:parent}
	Given a tree $\mathcal{G}_t=(\mathcal{V}_t,\mathcal{E}_t)$ with root vertex $r\in\mathcal{V}_t$, vertex $v\in\mathcal{V}_t$ is a parent of vertex $v'\in\mathcal{V}_t$ if and only if $\mathcal{L}_{v'}=\mathcal{L}_v+1$, where $\mathcal{L}_{v'}$ and $\mathcal{L}_v$ are the levels of $v'$ and $v$ in $\mathcal{G}_t$, respectively.
\end{definition}

\begin{definition}
	(Alternating tree) Given a matching $\mathcal{M}$, graph $\mathcal{G}$ is an alternating tree relative to $\mathcal{M}$ if and only if $\mathcal{G}$ is a tree with root vertex $r\in\mathcal{V}_t$, and for all vertices $v\in\mathcal{V}_t\backslash\{r\}$, the path between $v$ and $r$ is an alternating path relative to $\mathcal{M}$.
\end{definition}

Definition~\ref{def:parent} can be extended to the corresponding relationships of a child and a descendant of a vertex in an alternating tree. For example, vertex $v'\in\mathcal{V}_t$ is a child of vertex $v\in\mathcal{V}_t$ if and only if $v$ is a parent of $v'$. A vertex has only one parent but may have multiple children.

Consider a set of agents $\mathcal{A}=\{a_1,a_2,...,a_m\}$ and a set of tasks $\mathcal{B}:=\{b_1,b_2...,b_n\}$, with $m,n\in \mathbb{N}$. Assume that $\mathcal{A}\cap \mathcal{B}=\emptyset$ and that $m\geq n$. We represent all possible assignments of tasks to agents as a graph, i.e., let $\mathcal{G}_b=(\mathcal{V}_b,\mathcal{E}_b)$ be a complete weighted bipartite graph with vertex set $\mathcal{V}_b=\mathcal{A}\cup \mathcal{B}$ and edge set $\mathcal{E}_b=\{\{i,j\}|i\in\mathcal{A}, j\in\mathcal{B}\}$. Each edge $\{i,j\}\in \mathcal{E}_b$ is associated with a weight $\tau_{ij}\in\mathbb{R}$, which corresponds to the cost for agent $i$ carrying out task $j$. We define the set of all edge weights $\mathcal{W}=\{\tau_{ij}\in\mathbb{R}|i\in\mathcal{A},j\in \mathcal{B}\}$ and the function $w:\mathcal{E}_b\rightarrow \mathbb{R}$ that maps the edges of $\mathcal{G}_b$ to their weights. Let $\mathcal{C}(\mathcal{G}_b)$ be the set of all MCMs of $\mathcal{G}_b$. The BAP for graph $\mathcal{G}_b$ is formulated as
\begin{align} \label{eq:BAP}
\text{BAP}:\qquad&\min_{\mathcal{M}\in \mathcal{C}(\mathcal{G}_b)}\max_{e\in \mathcal{M}} \quad w(e).
\end{align}

We also define the set of MCMs that are solutions to the BAP given in~(\ref{eq:BAP}).

\begin{definition}[Bottleneck assignment] \label{def:botass}
	The set of all bottleneck assignments of $\mathcal{G}_b$ is $\mathcal{S}(\mathcal{G}_b):=\arg\min_{\mathcal{M}\in \mathcal{C}(\mathcal{G}_b)}\max_{e\in \mathcal{M}}w(e)$.
\end{definition}

\begin{definition}[Bottleneck edge] \label{def:botedge}
	Given any bottleneck assignmnent $\mathcal{M}\in\mathcal{S}(\mathcal{G}_b)$, any $e\in \arg\max_{e'\in\mathcal{M}} w(e')$ is a bottleneck edge of $\mathcal{G}_b$.
\end{definition}

For the forthcoming analyses on approaches to solve the BAP and on partitioning agents, we introduce the novel concepts of a pruned edge set, a bottleneck cluster, a critical bottleneck edge, and task and agent trees. Definition~\ref{def:phi} introduces a pruned edge set containing all edges in $\mathcal{M}$ and all edges that have weight strictly smaller than the largest edge in $\mathcal{M}$. With this tool we then define a bottleneck cluster and a critical bottleneck edge.

\begin{definition}[Pruned edge set] \label{def:phi}
	Given an MCM $\mathcal{M}$ of graph $\mathcal{G}_b$, a pruned edge set is defined as $\phi(\mathcal{G}_b,\mathcal{M}):=\mathcal{M}\cup \{e\in\mathcal{E}_b|w(e)<\max_{e'\in\mathcal{M}} w(e')\}$.
\end{definition}

There is existing literature on clustering agents and forming teams of agents that coordinate to carry out complex tasks, e.g.,~\cite{team1,team2}. A bottleneck cluster, defined below, is used to derive conditions for efficiently solving two separate BAPs as a combined problem. In the context of team formation, it provides a novel way to form a team of agents based on knowledge of the bottleneck edge.

\begin{definition}[Bottleneck cluster] \label{def:cluster}
	Consider a bottleneck assignment $\mathcal{M}\in\mathcal{S}(\mathcal{G}_b)$ and a bottleneck edge $e\in\mathcal{M}$ of $\mathcal{G}_b$. Graph $\mathcal{G}_b$ is a bottleneck cluster relative to $e$ if and only if for any vertices $v\in\mathcal{V}_b$, there exists an alternating path $\mathcal{P}$ relative to $\mathcal{M}$ between $v$ and a vertex $v'\in e$ such that $\mathcal{P}\subseteq\phi(\mathcal{G}_b,\mathcal{M})$.
\end{definition}

\begin{definition}[Critical bottleneck edge] \label{def:crit}
	Given an MCM $\mathcal{M}$ of graph $\mathcal{G}_b$, edge $e$ is a critical bottleneck edge of $\mathcal{G}_b$ relative to $\mathcal{M}$ if and only if $e\in\arg\max_{e\in\mathcal{M}} w(e)$ and $\phi(\mathcal{G}_b,\mathcal{M})\backslash\{e\}$ does not contain an augmenting path relative to $\mathcal{M}\backslash \{e\}$.
\end{definition}



In the following section, we show that all critical bottleneck edges are bottleneck edges. We can construct two particular subgraphs of a bottleneck cluster by partitioning its vertices and edges.

\begin{definition}[Agent and task trees] \label{def:tatree}
	Consider a bottleneck assignment $\mathcal{M}$ of graph $\mathcal{G}_b$ and a critical bottleneck edge $e_c=\{a_c,b_c\}$ relative to $\mathcal{M}$. Let $\mathcal{G}_b$ be a bottleneck cluster with respect to $e_c$. The agent and task trees of $\mathcal{G}_b$ are defined as $\mathcal{T}_\mu(\mathcal{G}_b):=(\mathcal{V}_{\mu},\mathcal{E}_{\mu})$ and $\mathcal{T}_\nu(\mathcal{G}_{b}):=(\mathcal{V}_{\nu},\mathcal{E}_{\nu})$ respectively, where $\mathcal{T}_\mu$ and $\mathcal{T}_\nu$ are alternating trees with the following properties. The sets of edges satisfy $\mathcal{E}_\mu\cup\mathcal{E}_\nu = \phi(\mathcal{G}_b,\mathcal{M})$, $\mathcal{E}_\mu\cap\mathcal{E}_\nu=\emptyset$, the sets of vertices satisfy $\mathcal{V}_\mu\cup\mathcal{V}_\nu=\mathcal{V}_b$, $\mathcal{V}_\mu\cap\mathcal{V}_\nu=\emptyset$, $\mathcal{V}_{\mu}$ contains the bottleneck agent $a_c= e_c\cap\mathcal{A}$, and $\mathcal{V}_{\nu}$ contains the bottleneck task $b_c = e_c\cap\mathcal{B}$.
\end{definition}

Fig.~\ref{fig:trees} in Section~\ref{sec:structure} illustrates, amongst other concepts, the agent and task trees for a given graph $\mathcal{G}_b$ and bottleneck assignment $\mathcal{M}$.

\section{A Distributed Algorithm for Solving the BAP} \label{sec:dist}

The following assumptions model a distributed setting, where information available to each agent may be limited.

\begin{assumption} \label{as:distedge}
	Assume an agent $i\in\mathcal{A}$ has access to the set of incident edges in graph $\mathcal{G}_b$, $\mathcal{E}_i := \{\{i,j\}\in \mathcal{E}_b|j\in\mathcal{B}\}$, and the weight of each edge in $\mathcal{E}_i$, i.e., the tuple $\mathcal{W}_i:=(w(\{i,b_1\}),w(\{i,b_2\}),...,w(\{i,b_n\}))$.
\end{assumption}

Note that $\mathcal{E}_b=\bigcup_{i=1}^m \mathcal{E}_i$ and $\mathcal{E}_{v}\cap \mathcal{E}_{v'}=\emptyset$ for $v,v'\in\mathcal{A}$, $v\neq v'$.

\begin{assumption} \label{as:distcom}
	Let communication between agents be modelled by a time invariant, undirected and connected graph $\mathcal{G}_C=(\mathcal{A},\mathcal{E}_C)$, where the set of agents $\mathcal{A}$ is the vertex set, with edge set $\mathcal{E}_C$ and diameter $D$. Agents communicate synchronously, i.e., all agents share a global clock and at discrete time steps of the clock, agents exchange local information $\mathcal{I}_i\in \mathcal{E}_i\times \mathcal{W}_i$ with all their neighbours $i'\in N(\mathcal{G}_C,i)$.
\end{assumption}

Note that instead of local information, agents are also able to pass on information received from neighbours in a previous time step of the global clock.

\begin{remark} \label{rem:distmatch}
	An MCM $\mathcal{M}$ of $\mathcal{G}_b$ can be jointly stored by agents. To this end, for each agent $i\in\mathcal{A}$ that is not free according to Definition~\ref{def:free}, let agent $i$ store its matched task $m_i$ in $\mathcal{M}$, i.e., $\{i,m_i\}\in\mathcal{M}$. We henceforth use $m_{i'}=-1$ to denote the matched task of a free agent $i'\in\mathcal{A}$. The function $\mathcal{Q}:\mathcal{C}(\mathcal{G}_b)\times \mathcal{A}\rightarrow \mathcal{B}\cup\{-1\}$ maps agents to their matched tasks in an MCM of $\mathcal{G}_b$.
\end{remark}

The problem of solving the BAP in a distributed manner is formally stated as follows.

\begin{problem} \label{prob:distBAP}
	Under Assumptions~\ref{as:distedge} and \ref{as:distcom}, obtain a solution to the BAP given by~(\ref{eq:BAP}).
\end{problem}

\subsection{An Algorithm for Solving the BAP} \label{subsec:distBAP}

We first present Algorithm~\ref{alg:centralised}, which we refer to as pruneBAP, for solving the BAP. Then, we discuss properties of this algorithm, including how each subroutine of pruneBAP can be implemented in the distributed setting.

\subsubsection{Guaranteeing a Bottleneck Assignment}

To initialise pruneBAP, we require an arbitrary MCM $\mathcal{M}_0\in\mathcal{C}(\mathcal{G}_b)$ of $\mathcal{G}_b$, e.g., agents and tasks can be initially matched based on arbitrary indexing, $\mathcal{M}_0= \{\{a_p,b_q\}|a_p\in \mathcal{A}, b_q\in \mathcal{B},p=q\}$.

\begin{algorithm}[thpb]
	\caption{pruneBAP}
	Input: Graph $\mathcal{G}_b=(\mathcal{V}_b,\mathcal{E}_b)$, and an MCM $\mathcal{M}_0$.\\
	Output: An MCM $\mathcal{M}$ of $\mathcal{G}_b$ that is a minimiser of (1).
	\begin{algorithmic}[1]
		\State $\mathcal{M}\gets \mathcal{M}_0$
		\State $\texttt{matching\_exists}\gets \texttt{True}$
		\While {\texttt{matching\_exists}}
		\State $(\bar{e},w(\bar{e}))\gets \texttt{\textsc{MaxEdge}}(\mathcal{M})$ \label{line:max}
		\State $\bar{\mathcal{E}}\gets \phi(\mathcal{G}_b,\mathcal{M})\backslash \{\bar{e}\}$ \label{line:remedge}
		\State $\bar{\mathcal{M}}\gets \mathcal{M}\backslash \{\bar{e}\}$ \label{line:remmatch}
		\State $\mathcal{M}_\nu \gets \texttt{\textsc{AugPath}}(\bar{e},\bar{\mathcal{M}},(\mathcal{V}_b,\bar{\mathcal{E}}))$ \label{line:augpath}
		\If {$\mathcal{M}_\nu \neq \bar{\mathcal{M}}$} \Comment{Augmenting path exists} \label{line:mnu}
		\State $\mathcal{M}\gets \mathcal{M}_\nu$
		\Else
		\State $\texttt{matching\_exists} \gets \texttt{False}$
		\EndIf
		\EndWhile
		\State \Return $\mathcal{M}$
	\end{algorithmic} \label{alg:centralised}
\end{algorithm}

Proposition~\ref{lem:distBAP} relates the existance of an augmenting path to the existance of an MCM within a set of edges and is crucial for proving that pruneBAP returns a bottleneck assignment in Proposition~\ref{lem:algcorrect}.

\begin{proposition}[Proof in Appendix~\ref{ap:lem:distBAP}] \label{lem:distBAP}
	Consider an MCM $\mathcal{M}$ of graph $\mathcal{G}=(\mathcal{V},\mathcal{E})$, a set of edges $\mathcal{E}'$ such that $\mathcal{M}\subseteq\mathcal{E}'\subseteq\mathcal{E}$, and an edge $e\in \mathcal{M}$. An augmenting path $\mathcal{P}\subseteq \mathcal{E}' \backslash\{e\}$ exists relative to $\mathcal{M}\backslash\{e\}$ if and only if there exists an MCM $\mathcal{M}'$ of $\mathcal{G}$ such that $\mathcal{M}'\subseteq \mathcal{E}' \backslash\{e\}$.
\end{proposition}

From Definition~\ref{def:crit} and Proposition~\ref{lem:distBAP}, we have the following corollary. We provide a proof for this corollary as some steps are not self-evident.

\begin{corollary}[Proof in Appendix~\ref{ap:cor:crit}] \label{cor:crit}
	Every critical bottleneck edge of $\mathcal{G}_b$ is a bottleneck edge of $\mathcal{G}_b$.
\end{corollary}

In Section~\ref{subsec:maxedge} and Section~\ref{subsec:augpath}, we verify that the following two assumptions hold and the functions can be implemented in a distributed setting.

\begin{assumption} \label{as:maxedge}
	Given an MCM $\mathcal{M}$ of $\mathcal{G}_b$, assume that there exists a function $\texttt{\textsc{MaxEdge}}(\mathcal{M})$ that returns the tuple $(\bar{e}, w(\bar{e}))$, where edge $\bar{e}\in\arg\max_{e\in\mathcal{M}} w(e)$.
\end{assumption}

\begin{assumption} \label{as:augpath}
	Given a graph $(\mathcal{V}_b,\mathcal{E}')$, an edge $\bar{e}\in\mathcal{E}'$, an MCM $\mathcal{M}$, and a matching $\bar{\mathcal{M}}=\mathcal{M}\backslash\{\bar{e}\}$, assume there exists function $\texttt{\textsc{AugPath}}(\bar{e},\bar{\mathcal{M}},(\mathcal{V}_b,\mathcal{E}'\backslash\{\bar{e}\}))$ that checks if there exists an augmenting path $\mathcal{P}$ relative to $\bar{\mathcal{M}}$ in $\mathcal{E}'\backslash\{\bar{e}\}$. If $\mathcal{P}$ exists, the function returns an MCM $\mathcal{M}_\nu=\bar{\mathcal{M}}\oplus \mathcal{P}$. If $\mathcal{P}$ does not exist, the function returns $\mathcal{M}_\nu=\bar{\mathcal{M}}$.
\end{assumption}

\begin{proposition}[Proof in Appendix~\ref{ap:lem:algcorrect}] \label{lem:algcorrect}
	If Assumptions~\ref{as:maxedge} and~\ref{as:augpath} hold, then pruneBAP returns a bottleneck assignment of $\mathcal{G}_b$.
\end{proposition}

In Fig.~\ref{fig:distBAP}, we show three iterations of the while-loop of pruneBAP for a toy BAP.

\begin{figure}[thpb]
	\centering
	\begin{tikzpicture}[scale=0.55]
	\tikzstyle{every node}=[font=\small]
	\node [] at (1,1) {Iteration 1:};
	\node [] at (0,0) (w1) {$e_{24}$};
	\node [] at (1,0) (w2) {$e_{42}$};
	\node [] at (2,0) (w3) {$e_{43}$};
	\node [] at (3,0) (w4) {$e_{34}$};
	\node [] at (4,0) (w5) {$e_{12}$};
	\node [] at (5,0) (w6) {$e_{21}$};
	\node [] at (6,0) (w7) {$e_{13}$};
	\node [] at (7,0) (w8) {$e_{22}$};
	\node [] at (8,0) (w9) {$e_{33}$};
	\node [] at (9,0) (w10) {$e_{23}$};
	\node [] at (10,0) (w11) {$e_{14}$};
	\node [] at (11,0) (w12) {$e_{31}$};
	\node [] at (12,0) (w13) {$e_{11}$};
	\node [] at (13,0) (w14) {$e_{41}$};
	\node [] at (14,0) (w15) {$e_{32}$};
	\node [] at (15,0) {$e_{44}$};
	
	\node [draw,circle,minimum size = 0.5cm] at (15,0) {};
	\node [draw,circle,minimum size = 0.5cm] at (12,0) {};
	\node [draw,circle,minimum size = 0.5cm] at (8,0) {};
	\node [draw,circle,minimum size = 0.5cm] at (7,0) {};
	\draw [dashed] (15.5,0.5) -- (15.5,-0.5);
	\node [] at (1,-1) {Iteration 2:};
	\node [] at (0,-2) (w1) {$e_{24}$};
	\node [] at (1,-2) (w2) {$e_{42}$};
	\node [] at (2,-2) (w3) {$e_{43}$};
	\node [] at (3,-2) (w4) {$e_{34}$};
	\node [] at (4,-2) (w5) {$e_{12}$};
	\node [] at (5,-2) (w6) {$e_{21}$};
	\node [] at (6,-2) (w7) {$e_{13}$};
	\node [] at (7,-2) (w8) {$e_{22}$};
	\node [] at (8,-2) (w9) {$e_{33}$};
	\node [] at (9,-2) (w10) {$e_{23}$};
	\node [] at (10,-2) (w11) {$e_{14}$};
	\node [] at (11,-2) (w12) {$e_{31}$};
	\node [] at (12,-2) (w13) {$e_{11}$};
	\node [] at (13,-2) (w14) {$e_{41}$};
	\node [] at (14,-2) (w15) {$e_{32}$};
	\node [] at (15,-2) {$e_{44}$};
	
	\node [draw,circle,minimum size = 0.5cm] at (12,-2) {};
	\node [draw,circle,minimum size = 0.5cm] at (0,-2) {};
	\node [draw,circle,minimum size = 0.5cm] at (8,-2) {};
	\node [draw,circle,minimum size = 0.5cm] at (1,-2) {};
	\draw [dashed] (12.5,-1.5) -- (12.5,-2.5);
	
	\node [] at (1,-3) {Iteration 3:};
	\node [] at (0,-4) (w1) {$e_{24}$};
	\node [] at (1,-4) (w2) {$e_{42}$};
	\node [] at (2,-4) (w3) {$e_{43}$};
	\node [] at (3,-4) (w4) {$e_{34}$};
	\node [] at (4,-4) (w5) {$e_{12}$};
	\node [] at (5,-4) (w6) {$e_{21}$};
	\node [] at (6,-4) (w7) {$e_{13}$};
	\node [] at (7,-4) (w8) {$e_{22}$};
	\node [] at (8,-4) (w9) {$e_{33}$};
	\node [] at (9,-4) (w10) {$e_{23}$};
	\node [] at (10,-4) (w11) {$e_{14}$};
	\node [] at (11,-4) (w12) {$e_{31}$};
	\node [] at (12,-4) (w13) {$e_{11}$};
	\node [] at (13,-4) (w14) {$e_{41}$};
	\node [] at (14,-4) (w15) {$e_{32}$};
	\node [] at (15,-4) {$e_{44}$};
	
	\node [draw,circle,minimum size = 0.5cm] at (4,-4) {};
	\node [draw,circle,minimum size = 0.5cm] at (5,-4) {};
	\node [draw,circle,minimum size = 0.5cm] at (3,-4) {};
	\node [draw,circle,minimum size = 0.5cm] at (2,-4) {};
	\draw [dashed] (5.5,-3.5) -- (5.5,-4.5);
	\end{tikzpicture}
	\caption{A demonstration of pruneBAP with agents $\mathcal{A}=\{a_1,a_2,a_3,a_4\}$ and tasks $\mathcal{B}=\{b_1,b_2,b_3,b_4\}$. The set of edges $\mathcal{E}_b$ is arranged in order of ascending weight, where $e_{pq}$ is the edge between agent $a_p$ and task $b_q$. At iteration 1, an initial arbitrary MCM is chosen, denoted by the 4 circled edges. In each iterations, edges to the right of the dashed lines are pruned from $\mathcal{E}_b$ and edges to the left form the set $\phi(\mathcal{G}_b,\mathcal{M})$. The weight of the largest edge in the MCMs is non-increasing, i.e., $w(e_{44})\geq w(e_{11})\geq w(e_{21})$. The algorithm terminates when it is not possible to form a matching of size 4 from the remaining edges to the left of the dashed line.}
	\label{fig:distBAP}
\end{figure}

\begin{remark} \label{rem:warmstart}
	Hueristically, the more edges we prune from the edge set $\mathcal{E}_b$, the faster pruneBAP converges to a solution, i.e., the fewer edges there are in $\phi(\mathcal{G}_b,\mathcal{M})$, the fewer possibilities there are to form an MCM of $\mathcal{G}_b$ from the edges within it. We apply this idea to the initialisation of pruneBAP. If pruneBAP is initialised with an arbitrary MCM $\mathcal{M}_0$, we call this as a cold-start. On the other hand, if some prior information can be exploited to choose a particular MCM with low edge weights, we call this initialisation a warm-start to pruneBAP.
\end{remark}

The warm-start property of pruneBAP is discussed in greater detail in the context of merging two BAPs in Section~\ref{sec:exploit}.

\subsection{Implementing pruneBAP with Distributed Computation}

We now discuss how the individual components of pruneBAP can be implemented with the distributed setting outlined in Assumptions~\ref{as:distedge} and~\ref{as:distcom} and Remark~\ref{rem:distmatch}. There are three main components to consider, i.e., finding the largest edge within a set of edges, removal of edges to form a pruned edge set, and searching for an augmenting path.

\subsubsection{Distributed Implementation of $\texttt{\textsc{MaxEdge}}()$} \label{subsec:maxedge}

We discuss a function that allows Assumption~\ref{as:maxedge} to be satisfied under Assumptions~\ref{as:distedge} and~\ref{as:distcom}. Given an edge set $\hat{\mathcal{E}}\subseteq \mathcal{E}_b$, agents must find an edge
\begin{equation} \label{eq:maxconsensus}
\bar{e}=\arg\max_{e\in\hat{\mathcal{E}}} w(e),
\end{equation}
without conflict. See Remark~\ref{rem:tiebreak} regarding conflict. The problem in~(\ref{eq:maxconsensus}) is known as the max-consensus problem \cite{max1,max2}.

\begin{lemma} [Proof in~\cite{max2}]\label{lem:maxconsensus}
	Under Assumptions~\ref{as:distedge} and~\ref{as:distcom}, the max-consensus algorithm in \cite[eq. (3)]{max2} solves the problem given by~(\ref{eq:maxconsensus}) and fulfills Assumption~\ref{as:maxedge}.
\end{lemma}

\begin{remark} \label{rem:tiebreak}
	Given an arbitrary set of edges $\hat{\mathcal{E}}\subseteq \mathcal{E}_b$, let $\alpha = \arg\max_{e\in\hat{\mathcal{E}}} w(e)$ be the set of all edges in $\hat{\mathcal{E}}$ with largest weight. If $\alpha$ is not a singleton, a deterministic method is required for selecting one edge to implement a function that satisfies Assumption~\ref{as:maxedge}. For instance, given the agents are indexed, the edge incident to the agent with the lower index number can be selected.
\end{remark}

Note that to satisfy Assumption~\ref{as:maxedge}, both the edge $\bar{e}$ and the corresponding weight $w(\bar{e})$ have to be returned. Therefore, agents $i\in\mathcal{A}$ must communicate the largest weight and the edge corresponding to it in each communication instance.

\subsubsection{Distributed Edge Removal}

Given the tuple $(\bar{e}, w(\bar{e}))$ from Assumption~\ref{as:maxedge} and given that each agent $i\in\mathcal{A}$ has access to $\mathcal{E}_i$ and $m_i$ from Assumption~\ref{as:distedge}, each agent $i$ can locally determine the set
\begin{equation} \label{eq:distphi}
\bar{\mathcal{E}}_i= \{i,m_i\}\cup\{e\in\mathcal{E}_i|w(e)<w(\bar{e})\}.
\end{equation}

Then, we have a distributed representation of the pruned edge set, since $\phi(\mathcal{G}_b,\mathcal{M})=\bigcup_{i\in\mathcal{A}}\bar{\mathcal{E}}_i$. To represent $\bar{\mathcal{E}}=\phi(\mathcal{G}_b,\mathcal{M})\backslash \{\bar{e}\}$ it remains for the particular agent $\bar{i}=\mathcal{A}\cap\bar{e}$ to additionally prune $\bar{e}=\{\bar{i},m_{\bar{i}}\}$, i.e.,
\begin{equation} \label{eq:distremedge}
\bar{\mathcal{E}}_{\bar{i}}= \{e\in\mathcal{E}_{\bar{i}}|w(e)<w(\bar{e})\}.
\end{equation}

Similarly edge $\bar{e}$ is locally removed from the matching $\mathcal{M}$ in Line~\ref{line:remmatch} of pruneBAP. To complete the distributed edge removal according to Remark~\ref{rem:distmatch}, agent $\bar{i}$ is labelled as free, i.e.,
\begin{equation} \label{eq:distremmatch}
m_{\bar{i}}=-1.
\end{equation}

\subsubsection{Distributed Implementation of $\texttt{\textsc{AugPath}}()$} \label{subsec:augpath}

We now present a distributed function satisfying Assumption~\ref{as:augpath}. Given an edge set $\bar{\mathcal{E}}$, the function searches for an augmenting path relative to a given matching $\bar{\mathcal{M}}$ under Assumptions~\ref{as:distedge} and~\ref{as:distcom}.

\makeatletter
\renewcommand{\ALG@name}{Function}
\makeatother
\addtocounter{algorithm}{-1}

\begin{algorithm}[thpb]
	\caption{Distributed DFS for finding an augmenting path.}
	Input: Edge $\{\bar{i},\bar{j}\}$, matching $\bar{\mathcal{M}}$, graph $(\mathcal{V}_b,\bar{\mathcal{E}})$.\\
	Output: New MCM $\mathcal{M}_\nu$.
	\begin{algorithmic}[1]
		\Function{AugDFS}{$\{\bar{i},\bar{j}\}$, $\bar{\mathcal{M}}$, $(\mathcal{V}_b,\bar{\mathcal{E}})$}
		\State $F\gets \emptyset$ \Comment{Set of explored agents}
		\State $m_i\gets \mathcal{Q}(\bar{\mathcal{M}},i),\forall i\in\mathcal{A}$ 
		\State $\nu_i\gets m_i,\forall i\in \mathcal{A}$ \label{line:nuinit}
		\State $\texttt{search\_complete}\gets \texttt{False}$
		\State $t\gets \bar{j}$
		\While {$\neg \texttt{search\_complete}$ }
		\State $S_t\gets \{i\in\mathcal{A}| i\notin F, \{i,t\}\in \bar{\mathcal{E}}\}$
		\State $a^*\gets a\in S_t$ \Comment{Choose next agent to explore} \label{line:chooseexplore}
		\If {$S_t=\emptyset$ and $t=\bar{j}$} \Comment{No remaining agents}
		\State $\texttt{search\_complete}\gets \texttt{True}$
		\ElsIf {$S_t=\emptyset$ and $t\neq \bar{j}$} \label{line:backtrack} \Comment{$t$ has no children}
		\State $a^*\gets k$, s.t. $m_k=t$ 
		\State $t\gets \nu_{a^*}$
		\State $\nu_{a^*}\gets m_{a^*}$ \label{line:nuback}		
		\ElsIf {$m_{a^*}=-1$} \Comment{Free agent found}
		\State $\nu_{a^*}\gets t$ \label{line:nudone}
		\State $\texttt{search\_complete}\gets \texttt{True}$
		\Else \label{line:forward} \Comment{Explore next agent}
		\State $\nu_{a^*}\gets t$ \label{line:nufor}
		\State $F\gets F\cup \{a^*\}$
		\State $t\gets m_{a^*}$ \label{line:childtask}
		\EndIf
		\EndWhile
		\State $\mathcal{M}_\nu\gets \{\{i,\nu_i\}|i\in\mathcal{A},\nu_i\neq -1\}$ \label{line:finish}
		\State \Return $\mathcal{M}_{\nu}$
		\EndFunction
	\end{algorithmic} \label{alg:augpath}
\end{algorithm}

In order to discuss Function~\ref{alg:augpath}, we introduce an alternating search that guarantees that all agents with an alternating path to the root are explored.

\begin{definition}[Alternating search] \label{as:tree}
	Given the complete bipartite graph $\mathcal{G}_b=(\mathcal{A}\cup\mathcal{B},\mathcal{E}_b)$ as introduced in Section~\ref{subsec:distBAP}, consider a subgraph $(\mathcal{A}\cup\mathcal{B},\bar{\mathcal{E}})$, where $\bar{\mathcal{E}}\subseteq\mathcal{E}_b$. Consider an MCM of $\mathcal{G}_b$, $\bar{\mathcal{M}}=\{\{a_1,m_{a_1}\}$, $\{a_2,m_{a_2}\}$, $...$, $\{a_n,m_{a_n}\}\}$, such that $\bar{\mathcal{M}}\backslash\{\{a_n,m_{a_n}\}\}\subseteq\bar{\mathcal{E}}$ and $\{a_n,m_{a_n}\}\notin\bar{\mathcal{E}}$. An alternating search is defined using the following vertex and edge set constructions.
	
	Consider a tree $(\mathcal{V}_1,\mathcal{E}_1)$ consisting of only the free root task vertex $m_{a_n}\in \mathcal{B}$, i.e., $\mathcal{V}_1=\{m_{a_n}\}$ and $\mathcal{E}_1=\emptyset$. For all $k\in\{1,2,...,f-1\}$, we construct vertex and edge sets $\mathcal{V}_{k+1}=\mathcal{V}_k\cup\{i_k, m_{i_k}\}$ and $\mathcal{E}_{k+1}=\mathcal{E}_k\cup \{\{i_k,j_k\},
	\{i_k,m_{i_k}\}\}$, satisfying the condition that agent $i_k\in \mathcal{A}\backslash \mathcal{V}_k$ and task $j_k\in \mathcal{B}\cap \mathcal{V}_k$ are neighbours in the graph $(\mathcal{A}\cup\mathcal{B},\bar{\mathcal{E}})$, i.e., $\{i_k,j_k\}\in \bar{\mathcal{E}}$. Let $f$ be the iteration at which for all tasks $j_f\in \mathcal{B}\cap \mathcal{V}_f$, there does not exist an agent $i_f\in \mathcal{A}\backslash \mathcal{V}_f$, where $\{i_f,j_f\}\in \bar{\mathcal{E}}$.
\end{definition}


\begin{proposition}[Proof in Appendix~\ref{ap:prop:systematic}] \label{prop:systematic}
	Given the tree $(\mathcal{V}_f,\mathcal{E}_f)$ constructed in Definition~\ref{as:tree}, the following statement holds. For any agent $i\in\mathcal{A}$, if there exists an alternating path between $i$ and $m_{a_n}$ containing only elements of $\bar{\mathcal{E}}$, then the agent is in the tree, i.e., $i\in\mathcal{V}_f$.
\end{proposition}

Function~\ref{alg:augpath} implements a depth-first search (DFS) that adheres to the pattern of agent exploration in Definition~\ref{as:tree}. See~\cite{planning} regarding DFS. The goal of the search is to find a free agent and the root $\bar{j}$ is a free vertex, so by construction the path between that free agent and the root is an augmenting path. On the other hand, Proposition~\ref{prop:systematic} guarantees that failure to find a free agent means an augmenting path does not exist.

\begin{remark} \label{rem:dfsstore}
	Rather than storing the entire alternating tree, Function~\ref{alg:augpath} only requires storage of the alternating path between the root and the current vertex $t$. The following is a suggested method to achieve this with distributed storage. All agents keep track of the current vertex $t$ by storing a First-In/Last-Out stack of tasks; a task is added to the stack when the search proceeds in Line~\ref{line:forward} and removed when the search backtracks in Line~\ref{line:backtrack}. To implement this, it is necessary for the explored agent $a^*$ to communicate the vertex $m_{a^*}$ to all agents.
\end{remark}

\begin{remark} \label{rem:greedypath}
	An approach for choosing the next agent to explore in Line~\ref{line:chooseexplore} in Function~\ref{alg:augpath} is to have $a^*\gets \arg\min_{i\in S_t} w(\{i,t\})$. This corresponds to the min-consensus problem, which is almost identical to the max-consensus problem in~(\ref{eq:maxconsensus}), only that the max is replaced with min. This approach greedily explores edges with smaller weights first.
\end{remark}

\begin{lemma}[Proof in Appendix~\ref{ap:lem:augpath}] \label{lem:augpath}
	Under Assumptions~\ref{as:distedge} and~\ref{as:distcom}, the function $\texttt{\textsc{AugDFS}}()$ in Function~\ref{alg:augpath} fulfills the conditions on $\texttt{\textsc{AugPath}}()$ and therefore guarantees that Assumption~\ref{as:augpath} is satisfied.
\end{lemma}

\subsubsection{Distributed Implementation of pruneBAP}
The algorithm pruneBAP returns a bottleneck assignment. Each component of pruneBAP can be implemented is a distributed setting provided by Assumptions~\ref{as:distedge} and~\ref{as:distcom}.

\begin{theorem}[Proof in Appendix~\ref{ap:thm:pruneBAP}] \label{thm:pruneBAP}
	The algorithm pruneBAP solves Problem~\ref{prob:distBAP}. 
\end{theorem}

Furthermore, we can bound the complexity of pruneBAP.

\begin{definition}[Time step] \label{def:timesteps}
	Given Assumption~\ref{as:distcom}, one time step refers to one time step or tick of the global clock shared by all agents.
\end{definition}

\begin{proposition}[Proof in Appendix~\ref{ap:prop:worst-case}] \label{prop:worst-case}
	The worst-case complexity, in terms of time steps, of the distributed implementation of pruneBAP is order $\mathcal{O}(mn^2D)$.
\end{proposition}

If $m=n$, then $|\mathcal{E}_b|=n^2$ and pruneBAP has a worst-case complexity of order $\mathcal{O}(n^3D)$. The number of iterations of the while-loop in Function~\ref{alg:augpath} is at most $2n-1$, where in the worst-case the search explores $n$ matched agents and backtracks $n-1$ times before terminating. The root $\bar{j}$ is the only free task vertex in the graph, which narrows down the search for an augmenting path. Without exploiting the fact that $\bar{j}$ known, finding an augmenting path in a bipartite graph in general has complexity $\mathcal{O}(n^2)$, from~\cite{assignprob,hk_alg}. We discuss the complexity of Function~\ref{alg:augpath} in more detail in the following section, where we make a comparison with an alternative method to implement the distributed search for an augmenting path. 

\section{Alternative Implementation of Distributed Augmenting Path Search} \label{sec:alternative}

Function~\ref{alg:augpath} conducts a distributed DFS that allows Assumption~\ref{as:augpath} to be satisfied. However, this is not the only function that fulfills the requirements of Assumption~\ref{as:augpath}. We introduce a distributed breadth-first search (BFS) approach in Function~\ref{alg:bfs} that also adheres to the construction of an alternating tree in Definition~\ref{as:tree}.

\subsection{BFS Augmenting Path Search}

Function~\ref{alg:bfs} does not implement a standard BFS~\cite{planning}. Multiple agents $i\in\mathcal{A}$ are explored per iteration of the while-loop starting on Line~\ref{line:while2} rather than only one agent at a time.

\begin{algorithm}[thpb]
	\caption{Distributed BFS for finding an augmenting path.}
	Input: Edge $\{\bar{i},\bar{j}\}$, matching $\bar{\mathcal{M}}$, graph $(\mathcal{V}_b,\bar{\mathcal{E}})$.\\
	Output: New MCM $\mathcal{M}_\nu$.
	\begin{algorithmic}[1]
		\Function{AugBFS}{$\{\bar{i},\bar{j}\}$, $\bar{\mathcal{M}}$, $(\mathcal{V}_b,\bar{\mathcal{E}})$}
		\State $F\gets \emptyset$ \Comment{Set of explored agents}
		\State $m_i\gets \mathcal{Q}(\bar{\mathcal{M}},i),\forall i\in\mathcal{A}$ 
		\State $\mathcal{R}\gets\emptyset$
		\State $\texttt{search\_complete}\gets \texttt{False}$
		\State $\bar{\mathcal{B}}\gets \{\bar{j}\}$ \Comment{Set of current tasks}
		\While {$\neg \texttt{search\_complete}$ } \label{line:while2}
		\State $\bar{S}\gets \{i\in\mathcal{A}| i\notin F, j\in \bar{\mathcal{B}},\{i,j\}\in \bar{\mathcal{E}}\}$ \label{line:barS}
		\State $F\gets F\cup \bar{S}$ \label{line:flag}
		\State $\nu_i\gets t\in \{j\in \bar{\mathcal{B}}|\{i,j\}\in\bar{\mathcal{E}}\},\forall i\in\bar{S}$ \label{line:parent}
		\State $\mathcal{R} \gets \mathcal{R} \cup \bigcup_{i\in \bar{S}} \{\{i,m_i\}, \{i,\nu_i\}\}$ \label{line:tree}
		\State $\bar{\mathcal{B}}\gets \bigcup_{i\in\bar{S}} m_i$ \label{line:tasks}
		\If {$\bar{S}=\emptyset$} \Comment{No remaining agents}
		\State $\texttt{search\_complete}\gets \texttt{True}$
		\State $\mathcal{M}_\nu\gets \bar{\mathcal{M}}$
		\ElsIf {$\exists i\in\bar{S},m_i=-1$} \label{line:foundfree} \Comment{Free agent found}
		\State $\texttt{search\_complete}\gets \texttt{True}$
		\State $a_{f}\in\{ i\in \bar{S}|m_i=-1\}$ \label{line:choosefree}
		\State $\mathcal{P}\gets P(\mathcal{R},a_f,\bar{j})$\label{line:paug}
		\State $\mathcal{M}_\nu\gets \bar{\mathcal{M}}\oplus \mathcal{P}$ \label{line:m_nu}
		\EndIf
		\EndWhile
		\State \Return $\mathcal{M}_{\nu}$
		\EndFunction
	\end{algorithmic} \label{alg:bfs}
\end{algorithm}

Each line in Function~\ref{alg:bfs} can be implemented with distributed computation, i.e., under Assumptions~\ref{as:distedge} and~\ref{as:distcom}. Function~\ref{alg:bfs} has two main steps. The first step is to determine the set of agents to explore next, i.e., to determine the elements of the set $\bar{S}$ in Line~\ref{line:barS}. To achieve this under a distributed setting, each unexplored agent individually checks for set membership. The second step is to store the explored agents, i.e, to store all the alternating paths in the alternating tree in Line~\ref{line:tree}. The following remark outlines a method to implement this step in the distributed setting.


\begin{remark} \label{rem:treestorage}
	Each time an agent $i\in \mathcal{A}$ is explored, agent $i$ stores its parent and child vertices, $\nu_i$ and $m_i$, respectively. This pair of tasks $(\nu_i,m_i)$ is then communicated to all other agents. All explored agents $i'\in F$ receive the pair of tasks and determine if task $m_i$ is to be added to their individual sets of descendents. Given that an agent $i'\in F$ receives a pair of tasks $(\nu_i,m_i)$, if $\nu_i$ is a descendant of $i'$, then $m_i$ is also a descendant of $i'$. When a free agent $a_f$ is explored, it only communicates its parent vertex $\nu_{a_f}$ to all other agents. The augmenting path between root $\bar{j}$ and $a_f$ is jointly stored by the corresponding agents that have $\nu_{a_f}$ in their set of descendants, i.e., for all agents $i\in F$, if $\nu_{a_f}$ is an element of the set of descendents of $i$, then $\mathcal{Q}(\mathcal{M}_\nu,i)=\nu_i$, where $\mathcal{Q}()$ is described in Remark~\ref{rem:distmatch}.
	
\end{remark}

\begin{lemma}[Proof in Appendix~\ref{ap:lem:BFS}] \label{lem:BFS}
	Under Assumptions~\ref{as:distedge} and~\ref{as:distcom}, the function $\texttt{\textsc{AugBFS}}()$ in Function~\ref{alg:bfs} fulfills the requirements for $\texttt{\textsc{AugPath}}()$ and therefore guarantees that Assumption~\ref{as:augpath} is satisfied.
\end{lemma}

Agents can execute Lines~\ref{line:barS} to~\ref{line:tasks} without waiting for other agents. The key point is that agents are explored in a way that adheres to Definition~\ref{as:tree}. This indicates that the synchronous communication between agents in Assumption~\ref{as:distcom} is a stronger assumption than necessary. However, we limit the discussion to the synchronous case, in this paper. 

\subsection{Comparing Distributed Search Methods}

TABLE~\ref{fig:table} shows a comparison of $\texttt{AugDFS}()$ given in Function~\ref{alg:augpath} to $\texttt{AugBFS}()$ given in Function~\ref{alg:bfs}. For the following discussion, iterations refers to the number of iterations of the while-loops of Function~\ref{alg:augpath} or~\ref{alg:bfs}, $n$ is the cardinality of an MCM of the given bipartite graph, and $D$ is the diameter of the communication graph.

\begin{table}[thpb]
	\centering
	\begin{tabular}{ p{2.57cm} || p{3.0cm} | p{2.1cm} } 
		Function & $\texttt{AugDFS}()$ & $\texttt{AugBFS}()$\\
		\hline\hline
		Agents explored per iteration & 1 (can be greedily chosen) & Multiple\\
		\hline
		Iterations to terminate & $2n-1$ & $n$\\
		\hline
		Time steps per iteration & $D$ & $D$\\
		\hline
		Information communicated per explored agent $i\in F$ & $m_i$, and $w(\{i,m_i\})$ for choosing the path greedily, and $i$ for tie-breaking & $m_i$, $\nu_i$, and $i$ for tie-breaking\\
		\hline
		Information stored by agent $i\in\mathcal{A}$ & $\nu_i$ and a First-In/Last-Out stack of tasks, see Remark~\ref{rem:dfsstore} & $\nu_i$ and the set of descendant tasks of $i$, see Remark~\ref{rem:treestorage}\\
	\end{tabular}
	\caption{Comparison of $\texttt{AugDFS}()$ and $\texttt{AugBFS}()$.} \label{fig:table}
\end{table}

In the function $\texttt{AugDFS}()$, only one vertex is explored in every iteration. To reach consensus on which agents to explore next, agents collect information from their neighbours, choose a local candidate and send their choice to their neighbours. In contrast, the function $\texttt{AugBFS}()$ explores all vertices with the same level in the tree at every iteration. Therefore, it returns the shortest augmenting path. However, since multiple agents are explored, agents collect information from their neighbours and then in turn send all this information to their neighbours. Thus, there is a trade-off between exploring multiple vertices per iteration and communicating more information between agents per iteration.

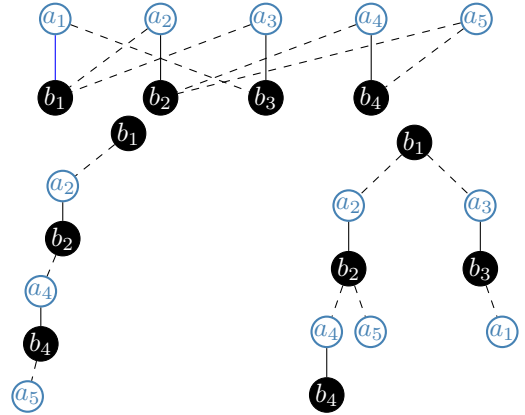
\begin{figure}[thpb]
	\centering
	\begin{tikzpicture}[scale=0.7, auto,swap]
	\node[agent] (a1) at (0,0) {$a_1$};
	\node[agent] (a2) at (2,0) {$a_2$};
	\node[agent] (a3) at (4,0) {$a_3$};
	\node[agent] (a4) at (6,0) {$a_4$};
	\node[agent] (a5) at (8,0) {$a_5$};
	
	\node[task] (b1) at (0,-1.5) {$b_1$};
	\node[task] (b2) at (2,-1.5) {$b_2$};
	\node[task] (b3) at (4,-1.5) {$b_3$};
	\node[task] (b4) at (6,-1.5) {$b_4$};
	
	\draw [blue] (a1) edge node[font=\tiny,left] {} (b1);
	\draw (a2) edge node[font=\tiny,left] {} (b2);
	\draw (a3) edge node[font=\tiny,left] {} (b3);
	\draw (a4) edge node[font=\tiny,left] {} (b4);
	\draw [dashed] (a2) edge (b1) (a3) edge (b1);
	\draw [dashed] (a4) edge node[font=\tiny,left] {} (b2);
	\draw [dashed] (a1) edge (b3) (a5) edge (b2);
	\draw [dashed] (a5) edge (b4);
	\end{tikzpicture}\\
	\begin{tikzpicture}[-,>=stealth',level/.style={sibling distance = 2.5cm/#1,
		level distance = 1cm},scale=0.7]
	\node [task] at (0,0) (b1) {$b_1$}
	child{  node [agent] (a2) {$a_2$} edge from parent[dashed]
		child[solid]{ node [task] (b2) {$b_2$} edge from parent[solid]
			child{  node [agent] (a4) {$a_4$} edge from parent[dashed]
				child[solid]{ node [task] (b4) {$b_4$} edge from parent[solid]
					child{  node [agent] (a5) {$a_5$} edge from parent[dashed]
					}
					child{ node [] (a4) {} edge from parent[draw=none]}
				}
			}
			child{ node [] (a4) {} edge from parent[draw=none]}
		}
	}
	child{ node [] (a4) {} edge from parent[draw=none]}
	;
	\end{tikzpicture}\qquad\qquad
	\begin{tikzpicture}[-,>=stealth',level/.style={sibling distance = 2.5cm/#1,
		level distance = 1.2cm},scale=0.7]
	\node [task] at (0,0) (b1) {$b_1$}
	child{  node [agent] (a2) {$a_2$} edge from parent[dashed]
		child[solid]{ node [task] (b2) {$b_2$} edge from parent[solid]
			child{  node [agent] (a4) {$a_4$} edge from parent[dashed]
				child[solid]{ node [task] (b4) {$b_4$} edge from parent[solid]
				}
			}
			child{  node [agent] (a5) {$a_5$} edge from parent[dashed]
			}
		}
	}
	child{  node [agent] (a3) {$a_3$} edge from parent[dashed]
		child[solid]{ node [task] (b3) {$b_3$} edge from parent[solid]
			child{  node [] (asdfas) {} edge from parent[draw=none]
			}
			child{  node [agent] (a1) {$a_1$} edge from parent[dashed]
			}
		}
	}
	;
	\end{tikzpicture}
	\caption{Comparison of alternating trees constructed by the two different searches. Given an MCM $\mathcal{M}$ represented by the solid lines, the top bipartite graph shows the edges in a pruned edge set $\phi(\mathcal{G}_b,\mathcal{M})$. The edge highlighted in blue is the edge to be removed, $\bar{e}$, and $b_1$ is the root of the search for an augmenting path. The bottom left tree is constructed via $\texttt{AugDFS}()$ and the bottom right tree, via $\texttt{AugBFS}()$.} \label{fig:BFStree}
\end{figure}

Fig.~\ref{fig:BFStree} demonstrates the difference between the alternating trees constructed by the searches in $\texttt{AugDFS}()$ and $\texttt{AugBFS}()$. No backtracking is required in $\texttt{AugBFS}()$ but this may result in multiple solutions being found. We now consider the worst-case number of iterations for the searches to terminate. The worst-case for $\texttt{AugDFS}()$ is if all agents only have length-one paths to the root task vertex. Then, the search would explore $n$ times and backtrack $n-1$ times. On the other hand the worst-case for $\texttt{AugBFS}()$ is if graph $(\mathcal{V}_b,\bar{\mathcal{E}})$ is a path graph. Each agent has a different level, so the worst-case number of iterations is $n$. Whether or not an augmenting path exists, the worst-case number of iterations to conclude the search is lower for $\texttt{AugBFS}()$ than for $\texttt{AugDFS}()$.


%
%

\section{Solving the BAP by Combining Two Sub-BAPs} \label{sec:exploit}

We analyse structure of the BAP that allows efficient merging of two subproblems. Although in this paper we consider two subproblems, we note that the results can be generalised for merging multiple subproblems. We formalise the particular structure that is to be exploited with Problem~\ref{prob:sol1} below.

\begin{problem} \label{prob:sol1}
	Consider two sets of agents $\mathcal{A}_1=\{a_1,a_2,...,a_{m_1}\}$ and $\mathcal{A}_2=\{\alpha_1,\alpha_2,...,\alpha_{m_2}\}$ and two sets of tasks $\mathcal{B}_1=\{b_1,b_2,...,b_{n_1}\}$ and $\mathcal{B}_2=\{\beta_1,\beta_2,...,\beta_{n_2}\}$. Define the sets $\mathcal{A}_3:=\mathcal{A}_1\cup\mathcal{A}_2$ and $\mathcal{B}_3:=\mathcal{B}_1\cup\mathcal{B}_2$. Let $m_3=m_1+m_2$ and $n_3=n_1+n_2$ and assume $m_1\geq n_1$ and $m_2\geq n_2$. For all $k\in\{1,2,3\}$, define graph $\mathcal{G}_k:=(\mathcal{V}_k,\mathcal{E}_k)$, with vertex set $\mathcal{V}_k:=\mathcal{A}_k\cup \mathcal{B}_k$ and edge set $\mathcal{E}_k:=\{\{i,j\}|i\in\mathcal{A}_k,j\in\mathcal{B}_k\}$. The function $w:\mathcal{E}_3\rightarrow \mathbb{R}$ maps all edges to their weights. Given a bottleneck assignment $\mathcal{M}_i\in \mathcal{S}(\mathcal{G}_i)$ and bottleneck a edge $e_i\in\mathcal{M}_i$ of graph $\mathcal{G}_i$, for $i\in\{1,2\}$, find a bottleneck assignment of $\mathcal{G}_3$, i.e., find an MCM $\mathcal{M}_3\in \mathcal{S}(\mathcal{G}_3)$.
\end{problem}

%

\subsection{A Bound on the BAP Solution} \label{subsubsec:bound}

We derive an upper bound on the weight of a bottleneck edge of $\mathcal{G}_3$ that depends on the weights of the bottleneck edges $e_1$ and $e_2$ of $\mathcal{G}_1$ and $\mathcal{G}_2$.

\begin{lemma}[Proof in Appendix~\ref{ap:thm:bound}] \label{thm:bound}
	Consider the definitions given in Problem~\ref{prob:sol1}. Given any bottleneck edge $e_3$ of graph $\mathcal{G}_3$, the weight $w(e_3)\leq\max\{w(e_1),w(e_2)\}$.
\end{lemma}

The set $\mathcal{M}_1\cup\mathcal{M}_2$ is an MCM of $\mathcal{G}_3$. This MCM is not necessarily a bottleneck assignment of $\mathcal{G}_3$. However, the largest weight of any edge in $\mathcal{M}_1$ and $\mathcal{M}_2$ provides an upper bound on the bottleneck weight of $\mathcal{G}_3$ according to Lemma~\ref{thm:bound}. If this bound is low enough to be acceptable for a particular assignment appliction, there is no need to invest further computational resources to solve Problem~\ref{prob:sol1} exactly. If an exact solution is required, $\mathcal{M}_1\cup\mathcal{M}_2$ can be used to warm-start pruneBAP.


As noted in Remark~\ref{rem:warmstart}, warm-starting is a heuristic. A warm-start does not guarantee fewer iterations for convergence to a bottleneck assignment of $\mathcal{G}_3$. Fig.~\ref{fig:combine} demonstrates a warm-start to pruneBAP. The sets $\mathcal{A}_3$ and $\mathcal{B}_3$ are the same as the agent and task sets in Fig.~\ref{fig:distBAP}. In this example, a bottleneck assignment of $\mathcal{G}_3$ is found in 1 iteration of the while-loop of pruneBAP.
	
\begin{figure}[thpb]
	\centering
	\begin{tikzpicture}[scale=0.55]
	\tikzstyle{every node}=[font=\small]
	\node [] at (2.5,1) {Subproblem 1 with $\mathcal{M}_1$:};
	\node [] at (4,0) (w2) {$e_{12}$};
	\node [] at (5,0) (w3) {$e_{21}$};
	\node [] at (7,0) (w8) {$e_{22}$};
	\node [] at (12,0) (w9) {$e_{11}$};
	
	\node [draw,circle,minimum size = 0.5cm] at (4,0) {};
	\node [draw,circle,minimum size = 0.5cm] at (5,0) {};
	\draw [dashed] (5.5,0.5) -- (5.5,-0.5);
	\node [] at (2.5,-1) {Subproblem 2 with $\mathcal{M}_2$:};
	\node [] at (2,-2) (w1) {$e_{43}$};
	\node [] at (3,-2) (w4) {$e_{34}$};
	\node [] at (8,-2) (w6) {$e_{33}$};
	\node [] at (15,-2) (w7) {$e_{44}$};
	
	\node [draw,circle,minimum size = 0.5cm] at (2,-2) {};
	\node [draw,circle,minimum size = 0.5cm] at (3,-2) {};
	\draw [dashed] (3.5,-1.5) -- (3.5,-2.5);

	\node [] at (2.5,-3) {Warm-start, Iteration 1:};
	\node [] at (0,-4) (w1) {$e_{24}$};
	\node [] at (1,-4) (w2) {$e_{42}$};
	\node [] at (2,-4) (w3) {$e_{43}$};
	\node [] at (3,-4) (w4) {$e_{34}$};
	\node [] at (4,-4) (w5) {$e_{12}$};
	\node [] at (5,-4) (w6) {$e_{21}$};
	\node [] at (6,-4) (w7) {$e_{13}$};
	\node [] at (7,-4) (w8) {$e_{22}$};
	\node [] at (8,-4) (w9) {$e_{33}$};
	\node [] at (9,-4) (w10) {$e_{23}$};
	\node [] at (10,-4) (w11) {$e_{14}$};
	\node [] at (11,-4) (w12) {$e_{31}$};
	\node [] at (12,-4) (w13) {$e_{11}$};
	\node [] at (13,-4) (w14) {$e_{41}$};
	\node [] at (14,-4) (w15) {$e_{32}$};
	\node [] at (15,-4) {$e_{44}$};
	
	\node [draw,circle,minimum size = 0.5cm] at (2,-4) {};
	\node [draw,circle,minimum size = 0.5cm] at (3,-4) {};
	\node [draw,circle,minimum size = 0.5cm] at (4,-4) {};
	\node [draw,circle,minimum size = 0.5cm] at (5,-4) {};
	\draw [dashed] (5.5,-3.5) -- (5.5,-4.5);
	\end{tikzpicture}
	\caption{A demonstration of solving divided subproblems to warm-start a combined problem, with $\mathcal{A}_1=\{a_1,a_2\}$, $\mathcal{A}_2=\{a_3,a_4\}$, $\mathcal{B}_1=\{b_1,b_2\}$, and $\mathcal{B}_2=\{b_3,b_4\}$.}
	\label{fig:combine}
\end{figure}

\subsection{Conditions for Merging Two BAPs Efficiently} \label{sec:structure}

Before we describe how to merge two BAPs, we first derive properties of a bottleneck cluster, a critical bottleneck edge, and the agent and task trees of a bottleneck cluster. The following proposition shows how a critical bottleneck edge forms a particular alternating path between the bottleneck agent and bottleneck task.

\begin{proposition}[Proof in Appendix~\ref{ap:prop:cluster}] \label{prop:cluster}
	Consider a bottleneck assignment $\mathcal{M}$ of bipartite graph $\mathcal{G}_b$ and a critical bottleneck edge $e_c=\{a_c,b_c\}$ relative to $\mathcal{M}$. The path $\mathcal{P}=\{e_c\}$ is the unique alternating path in the pruned edge set $\phi(\mathcal{G}_b,\mathcal{M})$ relative to $\mathcal{M}$ between $a_c$ and $b_c$.
\end{proposition}

The following corollary describes the structure of a bottleneck cluster $\mathcal{G}_b$ with respect to $e_c$ based on Definitions~\ref{def:cluster} and~\ref{def:tatree} and Proposition~\ref{prop:cluster}. 

\begin{corollary} \label{cor:tree}
	Given a bottleneck assignment $\mathcal{M}$ of graph $\mathcal{G}_b$ and a critical bottleneck edge $e_c=\{a_c,b_c\}$ relative to $\mathcal{M}$, let $\mathcal{G}_b$ be a bottleneck cluster with respect to $e_c$ and let the agent and task trees of $\mathcal{G}_b$ be $\mathcal{T}_\mu(\mathcal{G}_b)=(\mathcal{V}_{\mu},\mathcal{E}_{\mu})$ and $\mathcal{T}_\nu(\mathcal{G}_{b})=(\mathcal{V}_{\nu},\mathcal{E}_{\nu})$. For all agents in the agent tree $a'\in\mathcal{V}_\mu\cap \mathcal{A}$ and for all tasks in the task tree $b'\in\mathcal{V}_\nu\cap \mathcal{B}$, we have that $\{a',b'\}\notin \phi(\mathcal{G}_b,\mathcal{M})$.
\end{corollary}

\begin{figure}[thpb]
	\centering
	\begin{tikzpicture}[-,>=stealth',level/.style={sibling distance = 2.5cm/#1,
		level distance = 1.2cm},scale=0.6]
	\node [task] at (0,0) (b1) {$b_1$}
	child{  node [agent] (a2) {$a_2$} edge from parent[dashed]
		child[solid]{ node [task] (b2) {$b_2$} edge from parent[solid]
			 node[right] {\small$3$}
		} node[left] {\small$2$}
	}
	child{ node [agent] (a3) {$a_3$} edge from parent[dashed]
		child[solid]{ node [task] (b3) {$b_3$} edge from parent[solid] node[right] {\small$7$}
		} node[right] {\small$5$}
	}
	child{ node [agent] (a4) {$a_4$} edge from parent[dashed]
		child[solid]{ node [task] (b4) {$b_4$} edge from parent[solid] node[left] {\small$20$}
		} node[right] {\small$11$}
	}
	;
	
	\node [agent] at (7,0) (a1) {$a_1$}
	child{  node [task] (b6) {$b_5$} edge from parent[dashed]
		child[solid]{ node [agent] (a6) {$a_5$} edge from parent[solid]
			node[right] {\small$20$}
		} node[left] {\small$13$}
	}
	child{ node [task] (b7) {$b_6$} edge from parent[dashed]
		child[solid]{ node [agent] (a7) {$a_6$} edge from parent[solid] node[left] {\small$19$}
		} node[left] {\small$17$}
	}
	child{ node [task] (b8) {$b_7$} edge from parent[dashed]
		child[solid]{ node [agent] (a8) {$a_7$} edge from parent[solid] node[left] {\small$5$}
		} node[right] {\small$5$}
	}
	;
	\path [draw] (b1) -- node [label=above:{\small$20$}] {} (a1);
	
	\end{tikzpicture}
	\caption{An example bottleneck cluster with respect to a critical bottleneck edge. Dotted lines represent edges not in the matching $\mathcal{M}$, solid lines represent edges in $\mathcal{M}$. A set of agents $\{a_1,a_2,...,a_7\}$ and a set of tasks $\{b_1,b_2,...,b_7\}$. Edge $\{a_1,b_1\}$ is the critical bottleneck edge described in Proposition~\ref{prop:cluster}.}
	\label{fig:trees}
\end{figure}
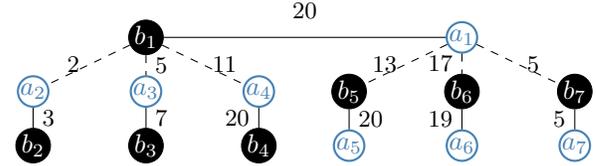

Fig.~\ref{fig:trees} illustrates Corollary~\ref{cor:tree}. The graph $\mathcal{G}_b$ is represented by two alternating trees $\mathcal{T}_\nu(\mathcal{G}_b)$ and $\mathcal{T}_\mu(\mathcal{G}_b)$ with roots $b_1$ and $a_1$, respectively. The vertices $a_1$ and $b_1$ are incident to the critical bottleneck edge $e_c$ with weight $w(e_c)=20$. Since $\mathcal{G}_b$ is a bottleneck cluster, all edges are elements of $\phi(\mathcal{G}_b,\mathcal{M})$.

In Section~\ref{sec:dist}, conditions for determining if an edge is a bottleneck edge of a given graph are proposed using Proposition~\ref{lem:distBAP}. We now build on this result and discuss corresponding conditions for determining if $\mathcal{M}_3=\mathcal{M}_1\cup\mathcal{M}_2$ is a bottleneck assignment of $\mathcal{G}_3$ when $\mathcal{G}_1$ and $\mathcal{G}_2$ are both bottleneck clusters.

\begin{lemma}[Proof in Appendix~\ref{ap:lem:reduce}] \label{lem:reduce}
	Consider the definitions given in Problem~\ref{prob:sol1}. Let $\mathcal{G}_1$ and $\mathcal{G}_2$ be bottleneck clusters with respect to $e_1$ and $e_2$, respectively. Let $e_1$ be a critical bottleneck edge of $\mathcal{G}_1$ relative to $\mathcal{M}_1$ and $e_2$ be a critical bottleneck edge of $\mathcal{G}_2$ relative to $\mathcal{M}_2$. Moreover assume $w(e_1)\geq w(e_2)$. If $w(e_3)<w(e_1)$, then the following conditions all hold:
	\begin{enumerate}[i.]
		\item there exists an edge in $\mathcal{E}_3$ with weight less than $w(e_1)$ between an agent $i\in\mathcal{A}_2$ and $b'\in\mathcal{V}_{\nu_1}\cap\mathcal{B}_1$, where $\mathcal{V}_{\nu_1}$ is the vertex set of the task tree $\mathcal{T}_\nu(\mathcal{G}_1)=(\mathcal{V}_{\nu_1},\mathcal{E}_{\nu_1})$ of $\mathcal{G}_1$;
		\item there exists an edge in $\mathcal{E}_3$ with weight less than $w(e_1)$ between a task $j\in\mathcal{B}_2$ and $a'\in\mathcal{V}_{\mu_1}\cap\mathcal{A}_1$, where $\mathcal{V}_{\mu_1}$ is the vertex set of the agent tree $\mathcal{T}_\mu(\mathcal{G}_1)=(\mathcal{V}_{\mu_1},\mathcal{E}_{\mu_1})$ of $\mathcal{G}_1$;
		\item there exists an alternating path $\mathcal{P}$ between $i$ and $j$ containing only edges with weight less than $w(e_1)$, and $|P\cap\mathcal{M}_2|>|P\backslash\mathcal{M}_2|$.
	\end{enumerate}
\end{lemma}

In general, the converse of Lemma~\ref{lem:reduce} does not hold when either $w(e_1)=w(e_2)$, or $\mathcal{M}_1$ contains more than one edge with weight equal to $w(e_1)$. We include these assumptions in the following theorem.

\begin{theorem}[Proof in Appendix~\ref{ap:thm:reduce2}] \label{thm:reduce2}
	Adopt the hypothesis of Lemma~\ref{lem:reduce}. Moreover, assume $w(e_1)> w(e_2)$ and $\arg\max_{\{i,j\}\in\mathcal{M}_1}w(\{i,j\})$ is a singleton. It holds that $w(e_3)<w(e_1)$ if and only if conditions i., ii., and iii. from Lemma~\ref{lem:reduce} hold.
\end{theorem}

\begin{figure}[thpb]
	\centering
	\begin{tikzpicture}
	\node[agent] at (-0.5,0) (a1) {$a_1$};
	\node[task] at (2,0) (b1) {$b_1$};
	
	\node[agent] at (-1,-2) (a2) {$a_2$};
	\node[task] at (-1,-1) (b2) {$b_2$};
	
	\node[agent] at (3,0.5) (a3) {$a_3$};
	\node[task] at (3,-0.5) (b3) {$b_3$};
	
	\node[agent] at (1,-3.5) (aa1) {$\alpha_1$};
	\node[task] at (0,-3) (bb1) {$\beta_1$};
	
	\node[agent] at (3,-2) (aa2) {$\alpha_2$};
	\node[task] at (2.2,-3) (bb2) {$\beta_2$};
	
	\path[draw,red] (a1) -- (b1);
	\path[draw] (a2) -- (b2);
	\path[draw] (a3) -- (b3);
	\path[draw] (aa1) -- (bb1);
	\path[draw] (aa2) -- (bb2);
	
	\path[draw,dashed] (a1) -- (b2);
	\path[draw,dashed] (b1) -- (a3);
	\path[draw,dashed] (aa1) -- (bb2);
	
	\path[draw,dashed,blue] (a2) -- (bb1);
	\path[draw,dashed,orange] (aa2) -- (b3);
	\end{tikzpicture}
	\caption{An illustration of Theorem~\ref{thm:reduce2}. The vertices are partitioned into two sets $\mathcal{V}_1=\{a_1,a_2,a_3\}\cup\{b_1,b_2,b_3\}$ and $\mathcal{V}_2=\{\alpha_1,\alpha_2\}\cup\{\beta_1,\beta_2\}$. Graph $\mathcal{G}_1$ is a bottleneck cluster with respect to $e_1=\{a_1,b_1\}$, and graph $\mathcal{G}_2$ is a bottleneck cluster with respect to $e_2=\{\alpha_2,\beta_2\}$. The length of each line corresponds to the weight of that edge. Solid lines show edges in $\mathcal{M}_1\cup\mathcal{M}_2$. The set of edges with dashed lines is an MCM with smaller weight than $w(e_1)$.}
	\label{fig:warmstart}
\end{figure}
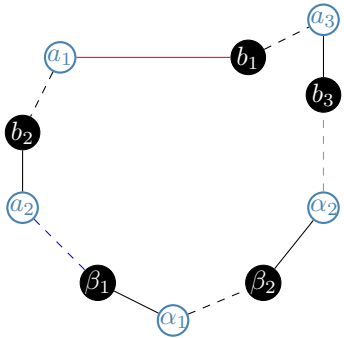

Fig.~\ref{fig:warmstart} illustrates a BAP where the conditions provided in Theorem~\ref{thm:reduce2} are satisfied. The orange dashed line is an edge that satisfies condition i.\ since $\alpha_2\in \mathcal{V}_2$ and there exists edge $\{\alpha_2,b_3\}\in \mathcal{V}_3$, where $b_3$ is a task in the vertex set of $\mathcal{T}_\nu(\mathcal{G}_1)$. The blue dashed line satisfies condition ii.\ since $\beta_1\in \mathcal{V}_2$ and there exists edge $\{\beta_1,a_2\}\in \mathcal{V}_3$, where $a_2$ is an agent in the vertex set of $\mathcal{T}_\mu(\mathcal{G}_1)$. Condition iii. is satisfied since there is an alternating path $\mathcal{P}=\{\{\alpha_1,\beta_1\},\{\alpha_2,\beta_2\}, \{\alpha_1,\beta_2\}\}$ between $\beta_1$ and $\alpha_2$, and $|\{\{\alpha_1,\beta_1\},\{\alpha_2,\beta_2\}\}|>|\{\{\alpha_1,\beta_2\}\}|$, i.e., $\mathcal{P}$ starts with a dashed line and ends with a dashed line.

Corollary~\ref{cor:main} follows from Lemma~\ref{thm:bound} and~\ref{lem:reduce}.

\begin{corollary} \label{cor:main}
	Adopt the hypothesis of Lemma~\ref{lem:reduce}. If at least one of conditions i., ii., or iii. in Lemma~\ref{lem:reduce} do not hold, then $\mathcal{M}_3=\mathcal{M}_1\cup\mathcal{M}_2$ is a solution to Problem~\ref{prob:sol1}. 
\end{corollary}

\section{Numerical Analysis} \label{sec:numerical}

In this section, we present a numerical analysis of the algorithm pruneBAP, the functions $\texttt{\textsc{AugDFS}}()$ and $\texttt{\textsc{AugBFS}}()$, and the warm-starting property. A time step is described in Definition~\ref{def:timesteps}. An iteration of pruneBAP represents the number of iterations of the while-loop of pruneBAP. Both appear as axis labels in several of the following plots.



\subsection{Numerical Analysis of pruneBAP}

We consider the following running scenario to illustrate the analyses in this subsection.

\subsubsection{Case Study}

Consider a equal number of agents and tasks, i.e.,  $n=m$, that are all represented by points in $\mathbb{R}^2$. Agents are to be assigned to move from their initial positions to assigned target positions such that a BAP with distance as weights is solved. To this end, we define the weights of the complete bipartite assignment graph $\mathcal{G}_b$ to be the Euclidean distance between agents and tasks. Unless otherwise specified, all coordinates  are generated from a uniform distribution between values of $0$ and $100$ normalised distance units.

By Assumption~\ref{as:distcom}, agents communicate synchronously and share a global clock. For the following examples, the communication between agents is modelled as a complete graph $\mathcal{G}_C$, i.e., all agents have a communication link to all other agents and the diameter $D=1$. For connected communication graphs that are not complete, the number of required time steps scales proportionally with the diameter $D$. Thus, the relative comparison of methods is fully illustrated by the case with fully connected communcation.



\subsubsection{Complexity of pruneBAP}

\begin{figure}
	\centering
	\includegraphics[scale=0.45]{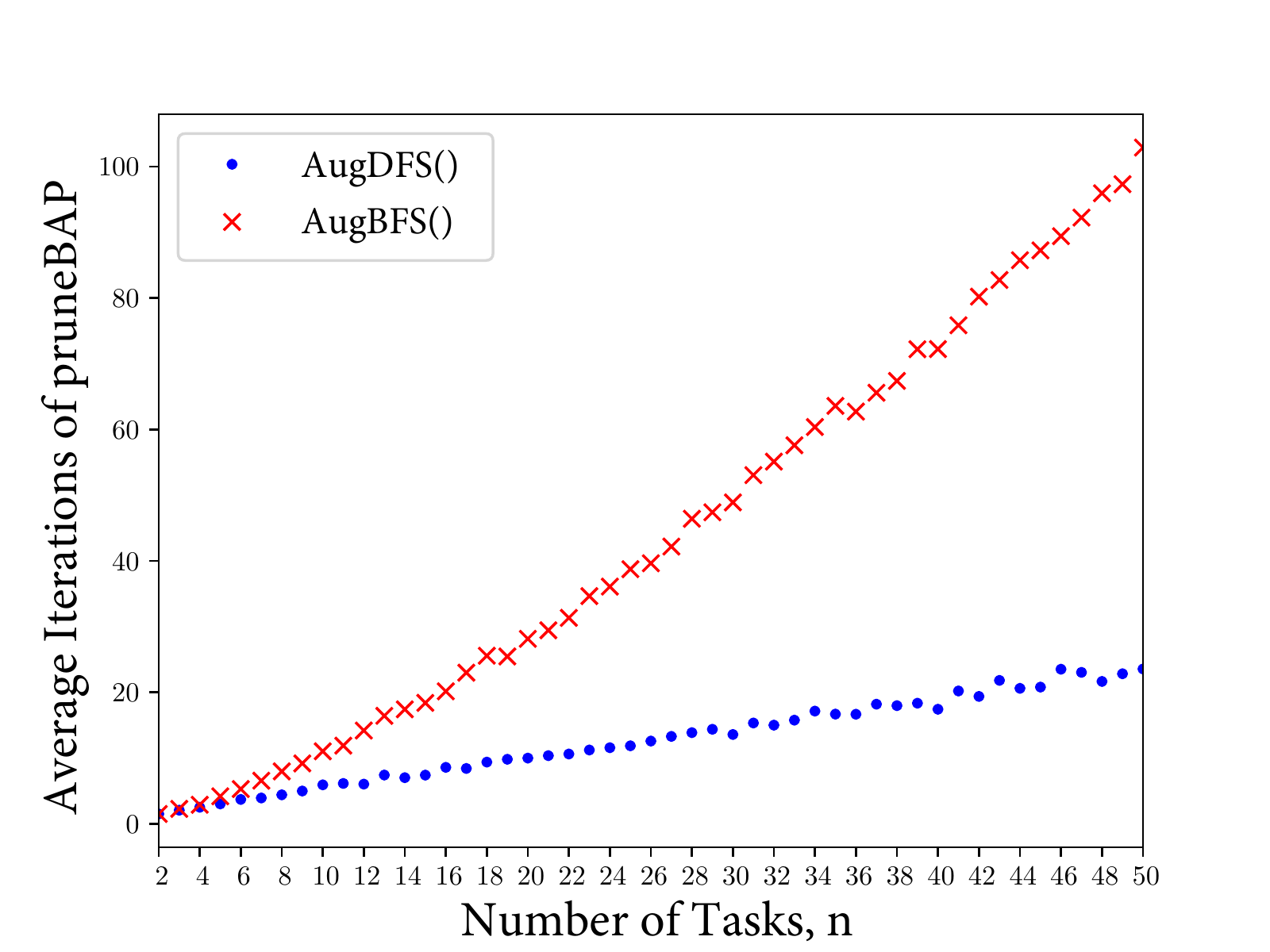}
	\caption{Average number of iterations of the while-loop of pruneBAP for completion versus the number of tasks, $n$.} \label{fig:itvsn}
\end{figure}

We evaluate the number of time steps it takes to run pruneBAP and compare the two implementations of an augmenting path search, $\texttt{\textsc{AugDFS}}()$ and $\texttt{\textsc{AugBFS}}()$. Fig.~\ref{fig:itvsn} shows the average number of iterations of the while-loop of pruneBAP for completion versus the number of tasks, $n$. For every value of $n$, the figure shows the average of $100$ realisations of the agent and task positions. As expected, $\texttt{\textsc{AugBFS}}()$ requires more iterations of the while-loop of pruneBAP than $\texttt{\textsc{AugDFS}}()$. This is because in $\texttt{\textsc{AugDFS}}()$ the augmenting path is constructed by a greedy minimisation of weights as described in Remark~\ref{rem:greedypath}. Therefore, the augmenting path that is found typically contains edges with small weights and fewer iterations of the while-loop of pruneBAP are required.

\begin{figure}[thpb]
	\centering
	\includegraphics[scale=0.45]{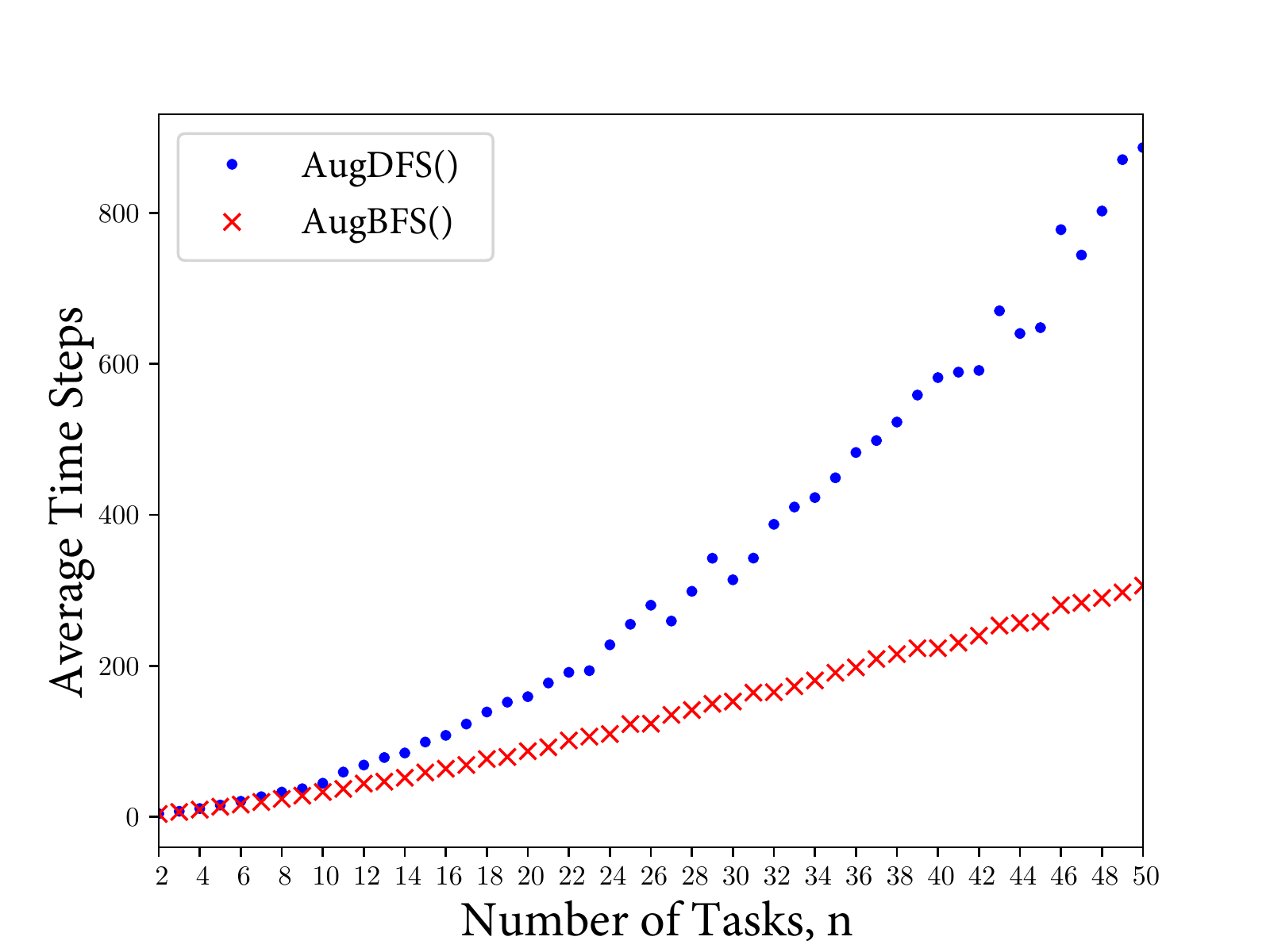}
	\caption{Average number of time steps for completion of pruneBAP versus the number of tasks, $n$.} \label{fig:stepsvsn}
\end{figure}

We recall that the worst-case complexity of pruneBAP is $\mathcal{O}(n^3)$. Fig.~\ref{fig:stepsvsn} shows the empirical average number of time steps for completion of pruneBAP versus the number of tasks $n$. Again, the values are averaged over 100 realisations of the agent and task positions. Although $\texttt{\textsc{AugDFS}}()$ results in fewer iterations of the while-loop of pruneBAP than $\texttt{\textsc{AugBFS}}()$, $\texttt{\textsc{AugBFS}}()$ results in fewer time steps for completion.

\subsubsection{Size of Messages Exchanged Between Agents}

TABLE~\ref{fig:table} states that each iteration of the while-loop of $\texttt{\textsc{AugDFS}}()$ or $\texttt{\textsc{AugBFS}}()$ consists of $D$ time steps, i.e., information from any one agent reaches all agents within $D$ time steps. Both functions require a similar amount of information to be exchanged per explored agent, so the number of explored agents per $D$ time steps is a proxy for the size of messages exchanged between agents.

\begin{figure}[thpb]
	\centering
	\includegraphics[scale=0.45]{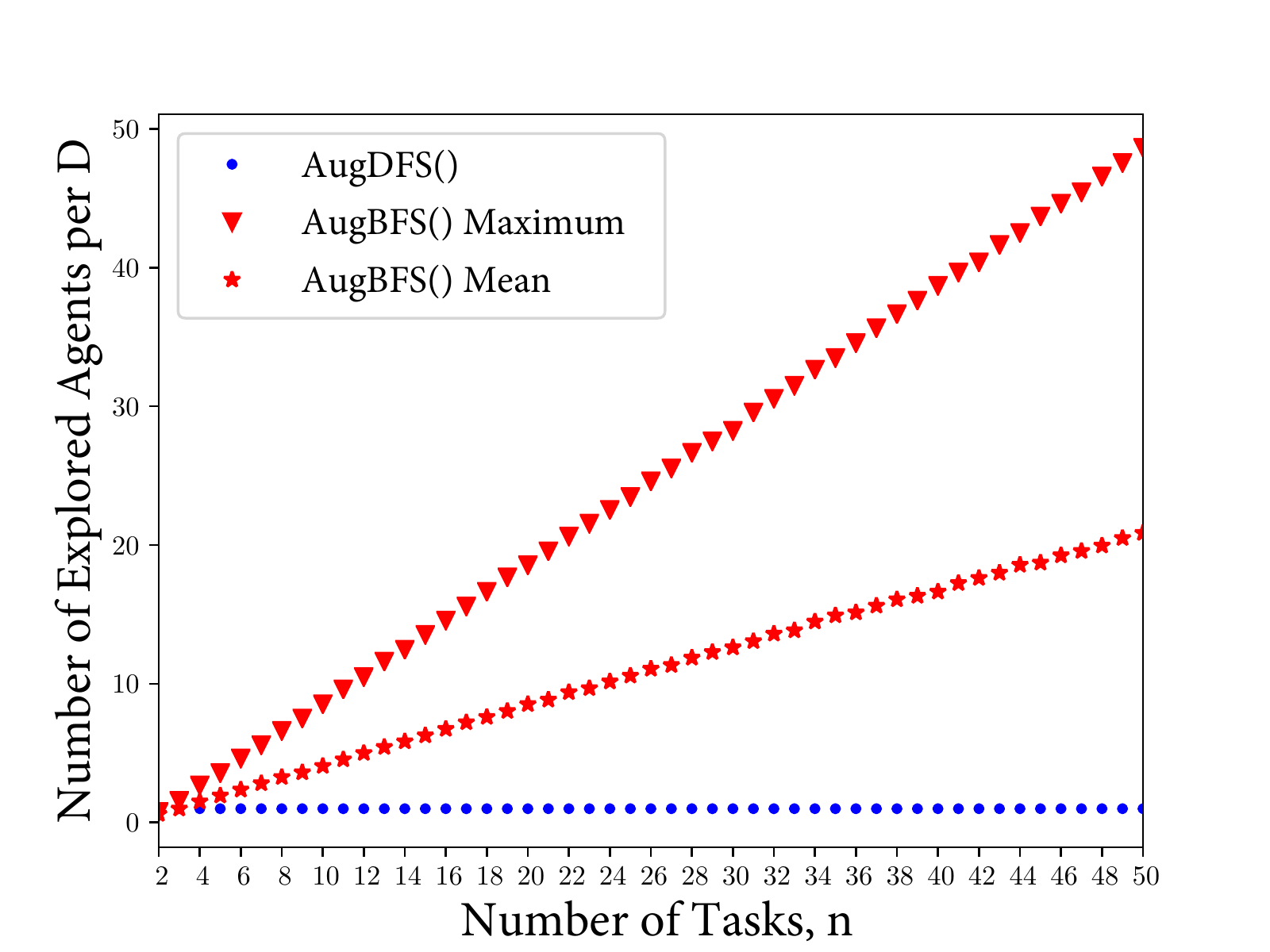}
	\caption{Average maximum number of explored agents and overall average number of explored agents per $D$ time steps for $\texttt{\textsc{AugBFS}}()$ versus number of tasks, $n$. For $\texttt{\textsc{AugDFS}}()$, one agent is explored per $D$ time steps.} \label{fig:mes}
\end{figure}

Note that $D=1$ is a special case where agents collect information from all agents and do not need to further relay any information. Fig.~\ref{fig:mes} shows the maximum number of explored agents and the mean number of explored agents per $D$ time steps for $\texttt{\textsc{AugBFS}}()$, as well as the exact number of explored agents per $D$ times steps for $\texttt{\textsc{AugDFS}}()$. Once again, the results are averaged over 100 realisations of the agent and task positions for every value of $n$.

In summary, pruneBAP converges faster when we apply $\texttt{\textsc{AugBFS}}()$ rather than $\texttt{\textsc{AugDFS}}()$, i.e., it requires fewer time steps on average to produce a solution to the BAP. However, the trade-off is that the messages exchanged between agents are larger and not of fixed length.




\subsubsection{Comparison Against Benchmark Distributed Assignment Method}

A greedy assignment algorithm is one that sequentially chooses the edge with lowest weight given prior selections while ensuring the chosen edges form a matching of $\mathcal{G}_b$. We compare the performance of pruneBAP with a distributed greedy algorithm. The greedy algorithm we consider as a benchmark is CBAA from~\cite{CBAA}, which can be applied in the distributed setting defined in Assumptions~\ref{as:distedge} and~\ref{as:distcom}.

\begin{figure}[thpb]
	\centering
	\includegraphics[scale=0.45]{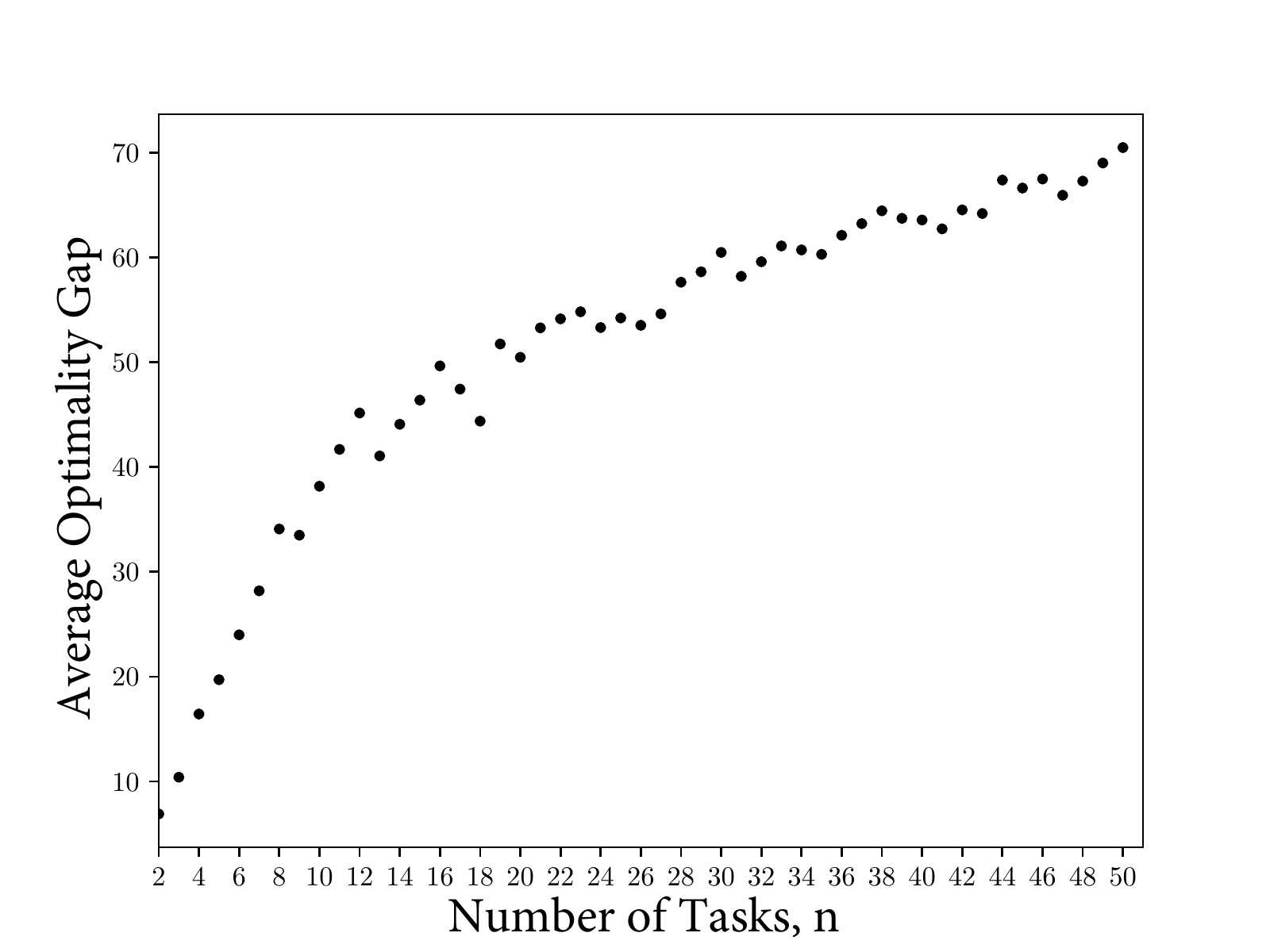}
	\caption{Average normalised optimality gap, $g-h$ versus $n$.} \label{fig:optimgap}
\end{figure}

A greedy approach does not necessarily solve the BAP. We first illustrate the optimality gap of CBAA. We denote the largest weight of all the edges corresponding to the assignment obtained from the CBAA approach as $g$ and the weight of a bottleneck edge as $h$. Fig.~\ref{fig:optimgap} shows the average optimality gap, $g-h$, averaged across 100 realisations of the agent and task positions per $n$. The optimality gap increases as $n$ increases. This happens because the weight of the largest edge in the final assignment found via pruneBAP decreases as $n$ increases. As $n$ increases, the number of MCMs in the problem increase, so it becomes more likely that an MCM with a smaller edge weights can be found. The weight $g$ for an assignment found via CBAA will not decrease as the number of tasks increases since it depends only on the last chosen edge, which has a expected weight value that is independent of $n$.

\begin{figure}[thpb]
	\centering
	\includegraphics[scale=0.45]{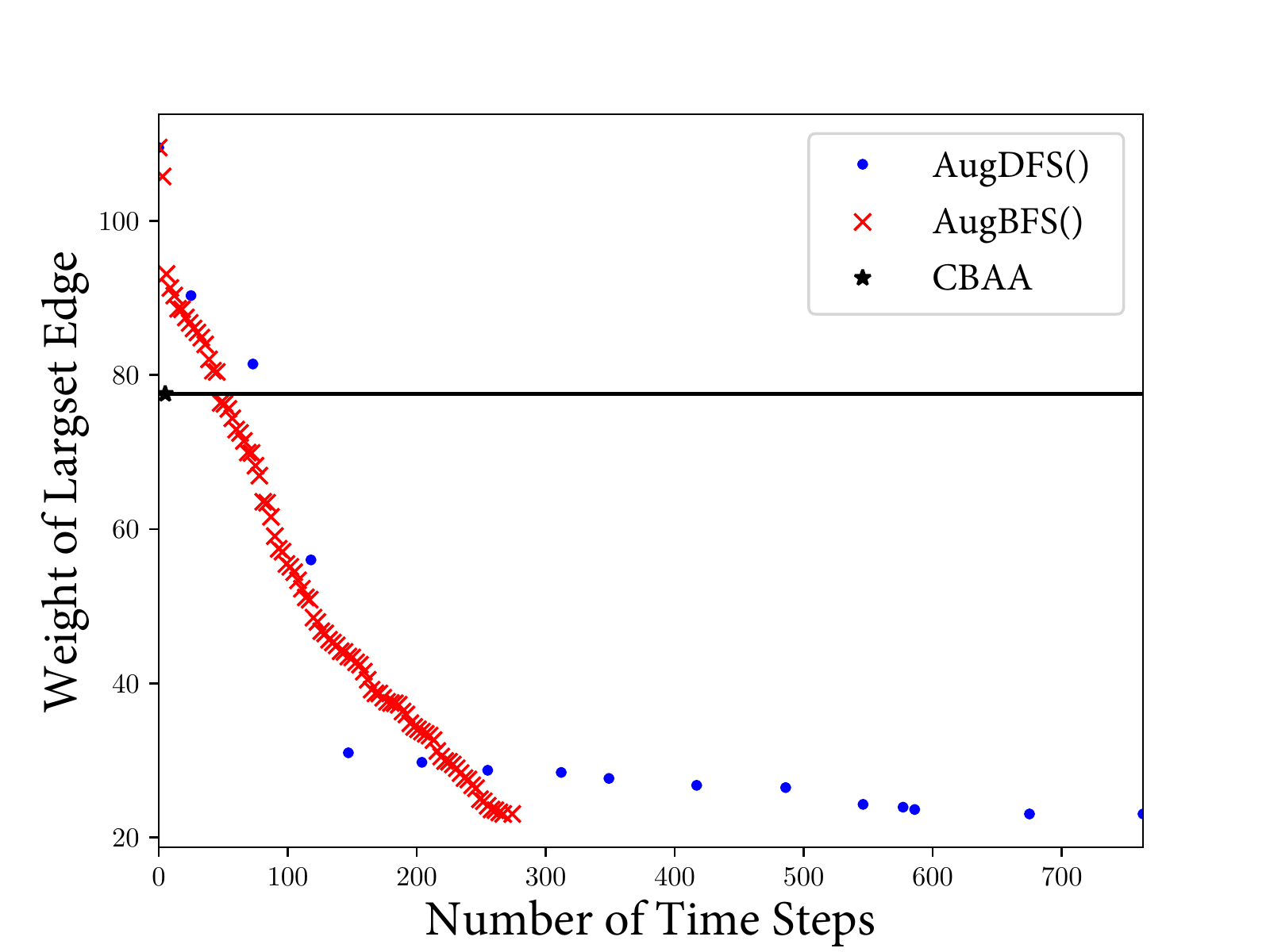}
	\caption{Weight of largest edge in MCM versus number of time steps.} \label{fig:kvsn}
\end{figure}

In Fig.~\ref{fig:kvsn}, we consider one example realisation of the agent and task positions with $m=n=50$. We investigate the number of time steps pruneBAP requires to find an MCM $\mathcal{M}$ with largest weight smaller than $g$. This is indicated by the time step at which the blue and red marks drop below the black line in Fig.~\ref{fig:kvsn}. Each point in Fig.~\ref{fig:kvsn} corresponds to an MCM. CBAA only produces an MCM in the final time step, whereas pruneBAP produces a series of MCMs with a non-increasing largest edge weight.

\begin{figure}[thpb]
	\centering
	\includegraphics[scale=0.45]{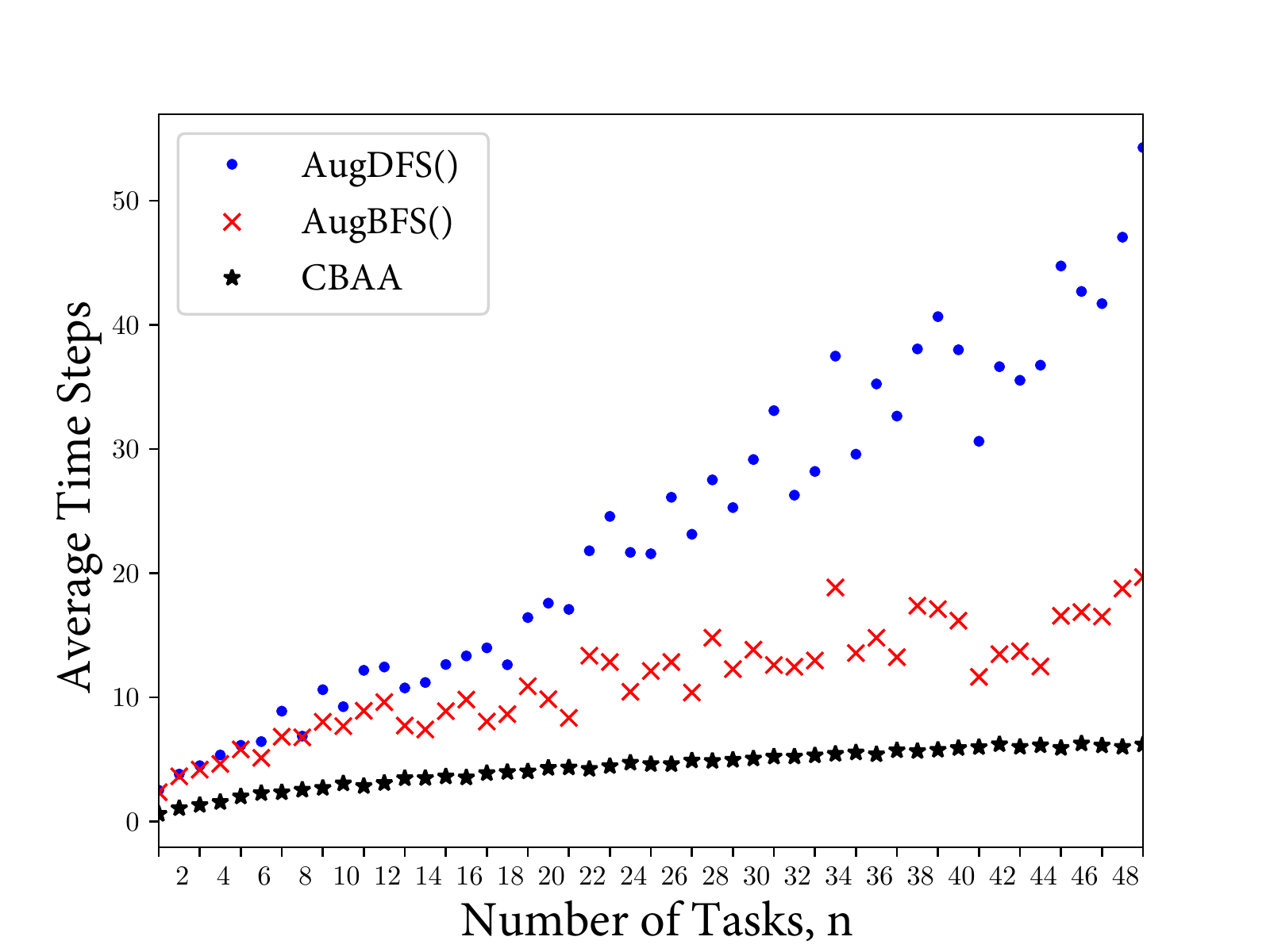}
	\caption{Average number of time steps for CBAA to converge and average number of time steps for pruneBAP to find an MCM with largest weight lower than the largest weight from the CBAA assignment versus the number of tasks $n$.} \label{fig:cvn}
\end{figure}

Fig.~\ref{fig:cvn} shows the average number of time steps required for pruneBAP to find an MCM $\mathcal{M}$ that has a largest weight smaller than $g$. Fig.~\ref{fig:cvn} also shows the number of time steps for CBAA to find an assignment. Once again, the averages were taken across 100 realisations for every $n$. 

Fig.~\ref{fig:cvn} suggests that on average, pruneBAP would converge faster by warm-starting with the assignment found via CBAA. Next, we consider another situation where warm-starting can be applied by exploiting the bottleneck assignments of subgraphs.

\subsubsection{Demonstration of Merging Subproblems Efficiently}
\begin{figure}[thpb]
	\centering
	\includegraphics[scale=0.45]{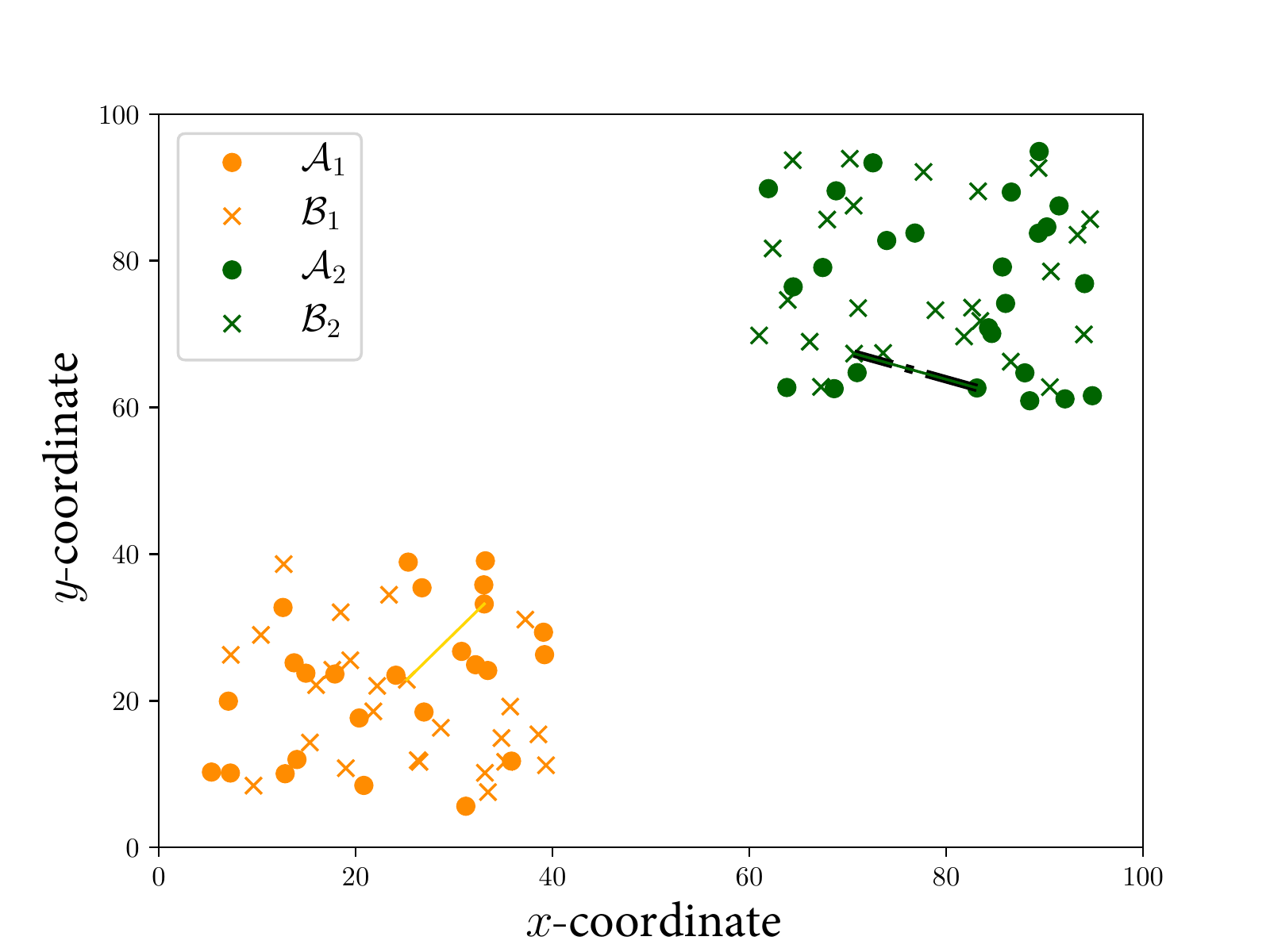}
	\caption{Case Study 2: Sample configuration of agents and tasks. The solid lines show the bottleneck edges from solving the BAP for graphs $\mathcal{G}_1$ and $\mathcal{G}_2$. The black dash-dotted line shows a bottleneck edge of $\mathcal{G}_3$.}
	\label{fig:cluster}
\end{figure}

Fig.~\ref{fig:cluster} shows two sets of agents, with initial locations given by the sets $\mathcal{A}_1=\{a_1,a_2,...,a_{25}\}$ and $\mathcal{A}_2=\{\alpha_1,\alpha_2,...,\alpha_{25}\}$ and two sets of goal locations given by the sets $\mathcal{B}_1=\{b_1,b_2,...,b_{25}\}$ and $\mathcal{B}_2=\{\beta_1,\beta_2,...,\beta_{25}\}$. Let the graphs $\mathcal{G}_1$, $\mathcal{G}_2$ and $\mathcal{G}_3$ be as defined in Problem~\ref{prob:sol1}. The coordinates are generated from uniform distributions between $5$ to $40$ and $60$ to $95$ normalised distance units. In this example, conditions i. and ii. in Lemma~\ref{lem:reduce} do not hold. Therefore, merging the bottleneck assignments of $\mathcal{G}_1$ and $\mathcal{G}_2$ results in a bottleneck assignment of $\mathcal{G}_3$.

Based on Fig.~\ref{fig:stepsvsn} and the worst-case complexity in Proposition~\ref{prop:worst-case}, more time steps are required for completion of pruneBAP for a larger number of tasks, $n$. Therefore, solving the two subproblems for graphs $\mathcal{G}_1$ and $\mathcal{G}_2$ takes fewer time steps than solving the full problem $\mathcal{G}_3$ from scratch, particularly when the subproblems are solved in parallel. The union of the MCMs of the two subproblems, i.e., $\mathcal{M}_1\cup\mathcal{M}_2$, can be used to warm-start pruneBAP to solve the combined problem with graph $\mathcal{G}_3$. In scenarios such as the one in shown in Fig.~\ref{fig:cluster}, where a condition in Lemma~\ref{lem:reduce} does not hold, warm-starting with $\mathcal{M}_1\cup\mathcal{M}_2$ leads pruneBAP to terminate after one iteration with $\mathcal{M}_3=\mathcal{M}_1\cup\mathcal{M}_2$.




\section{Conclusion}

We introduced tools, i.e., a pruned edge set, a bottleneck cluster, and a critical bottleneck edge to present an algorithm to solve the BAP. The pruneBAP algorithm iteratively produces an assignment with lower bottleneck weight with the final assignment being a solution to the BAP. The algorithm has several components, i.e., finding a pruned edge set, finding the largest edge amongst agents, and searching for an augmenting path. We derived methods to distribute the execution of each individual component over a network of agents with limited information. In particular, we compared two methods, $\texttt{\textsc{AugDFS}}()$ and $\texttt{\textsc{AugBFS}}()$, for conducting a distributed search for an augmenting path. The two methods provide a trade-off between computational complexity and amount of information that needs to be communicated between agents.

We identified structure that can be exploited when two BAPs are merged into a combined BAP. In that situation, the conditions on whether or not reassignment is necessary are related to the existance of an particular augmenting path between the two bottleneck clusters. We investigated the average numerical complexity of pruneBAP, $\texttt{\textsc{AugDFS}}()$, and $\texttt{\textsc{AugBFS}}()$. As a benchmark, we compared pruneBAP to CBAA, a greedy distributed assignment-finding algorithm. We provided an example where solving two separate BAPs leads to a solution of the combined problem.

Aside from their use for solving the BAP, augmenting path searches are applied in many approaches for solving other combinatorial optimisation problems. These include the problems of finding an MCM of a bipartite graph, finding the maximum cardinality intersection of a matroid, and finding the minimum of submodular functions, see~\cite{hk_alg,matroids,submodular,submodular2}. The distributed augmenting path search methods derived in this paper have the potential to provide benefits for solving other combinatorial optimisation problems. The extension of the methods to other applications is the subject of future work.

\bibliographystyle{unsrt}
\bibliography{mybib}  

\appendix
\section{Appendix}

\subsection{Auxiliary lemmas}
Lemmas~\ref{lem:aug} and~\ref{lem:berge} relate the existance of an augmenting path to the cardinality of matchings in a graph and are used to prove Proposition~\ref{lem:distBAP}. The notation $A\oplus B$ denotes the symmetric difference of the sets $A$ and $B$.

\begin{lemma}[Proof in \cite{hk_alg}] \label{lem:aug}
	Given a graph $\mathcal{G}$, if $\mathcal{M}$ is a matching of $\mathcal{G}$ and $\mathcal{P}$ is an augmenting path relative to $\mathcal{M}$, then $\mathcal{M}\oplus \mathcal{P}$ is also a matching of $\mathcal{G}$ and $|\mathcal{M}\oplus \mathcal{P}|=|\mathcal{M}|+1$.
\end{lemma}

\begin{lemma}[Berge's Theorem, Proof in \cite{hk_alg}] \label{lem:berge}
	Given a graph $\mathcal{G}$, a matching $\mathcal{M}$ is an MCM of $\mathcal{G}$ if and only if there is no augmenting path relative to $\mathcal{M}$ in the edge set of graph $\mathcal{G}$.
\end{lemma}

\subsection{Proof of Proposition~\ref{lem:distBAP}} \label{ap:lem:distBAP}
First, we prove sufficiency. Assume there exists an augmenting path $\mathcal{P}\subseteq \mathcal{E}'\backslash\{e\}$ relative to $\mathcal{M}\backslash\{e\}$. By Lemma~\ref{lem:aug}, $\mathcal{M}'=\mathcal{M}\backslash \{e\} \oplus \mathcal{P}$ is an MCM of $\mathcal{G}$. Since both $\mathcal{P}$ and $\mathcal{M}\backslash\{e\}$ are subsets of $\mathcal{E}'\backslash\{e\}$, their symmetric difference is also a subset of $\mathcal{E}'\backslash\{e\}$.

Next, we prove necessity. Assume there does not exist an augmenting path $\mathcal{P}\subseteq \mathcal{E}' \backslash\{e\}$ relative to $\mathcal{M}\backslash\{e\}$. By Lemma~\ref{lem:berge}, the matching $\mathcal{M}\backslash\{e\}$ is an MCM of the graph $(\mathcal{V},\bar{\mathcal{E}})$ with edge set $\bar{\mathcal{E}}= \mathcal{E}' \backslash\{e\}$. We know that $\mathcal{M}$ is an MCM of $\mathcal{G}$ and it has cardinality $|\mathcal{M}|$. Since $\mathcal{M}\backslash\{e\}$ has cardinality $|\mathcal{M}|-1$, it is not an MCM of $\mathcal{G}$. Thus, there does exist an MCM of $\mathcal{G}$ within the set $\mathcal{E}' \backslash\{e\}$.

\subsection{Proof of Corollary~\ref{cor:crit}} \label{ap:cor:crit}

Consider Proposition~\ref{lem:distBAP} with the graph $\mathcal{G}_b$, an MCM $\mathcal{M}$ of $\mathcal{G}_b$, the pruned edge set $\phi(\mathcal{G}_b,\mathcal{M})=\mathcal{E}'$, and a critical bottleneck edge $e$. Since $e$ is a critical bottleneck edge, there does not exist an MCM of $\mathcal{G}_b$ in $\phi(\mathcal{G}_b,\mathcal{M})\backslash\{e\}$, i.e., all MCMs of $\mathcal{G}_b$ must either contain $e$ or must contain edges with equal or larger weight that $w(e)$. Thus, $\mathcal{M}$ is a bottleneck assignment and $e$ is a bottleneck edge of $\mathcal{G}_b$.

\subsection{Proof of Proposition~\ref{lem:algcorrect}} \label{ap:lem:algcorrect}

Given Assumptions~\ref{as:maxedge} and~\ref{as:augpath}, we observe that Lines~\ref{line:remedge},~\ref{line:remmatch}, and~\ref{line:augpath} of pruneBAP test the sufficient and necessary conditions for the existance of the MCM $\mathcal{M}'$ based on Proposition~\ref{lem:distBAP}, with $e=\bar{e}$ and $\mathcal{E}' = \phi(\mathcal{G}_b,\mathcal{M})$.

Now, let $k\in\{1,2,...,f\}$ denote the iterations of the while-loop of Algorithm 1, where without loss of generality, $f$ is the final iteration. Let $\bar{\mathcal{E}}_k$ denote the set $\bar{\mathcal{E}}$ in Line~\ref{line:remedge} at iteration $k$ of the while-loop. By Proposition~\ref{lem:distBAP}, we have that for all iterations $k<f$, there exists an MCM $\mathcal{M}_k$ of $\mathcal{G}_b$ such that $\mathcal{M}_k\subseteq\bar{\mathcal{E}}_k$. At the final iteration $f$, an augmenting path is not found so pruneBAP returns $\mathcal{M}=\mathcal{M}_{f-1}$. It remains to show that the MCM $\mathcal{M}_{f-1}\subseteq\bar{\mathcal{E}}_{f-1}$ is a bottleneck assignment of $\mathcal{G}_b$. Assume for contradiction there exists an MCM $\mathcal{M}_f$ of $\mathcal{G}_b$ with all edges having weights strictly less than $e\in\arg\max_{e\in\mathcal{M}_{f-1}} w(e)$. This implies that there exists an MCM $\mathcal{M}_f$ of $\mathcal{G}_b$ such that $\mathcal{M}_{f}\subseteq \phi(\mathcal{G}_b,\mathcal{M}_{f-1})\backslash\{e\}=\bar{\mathcal{E}}_f$, which contradicts the assumption that $f$ is the final iteration. Following the above arguments, $e$ is a critical bottleneck edge of $\mathcal{G}_b$ and consequently $\mathcal{M}_{f-1}$ is a bottleneck assignment.

\subsection{Proof of Proposition~\ref{prop:systematic}} \label{ap:prop:systematic}

By construction, if task $j\in\mathcal{B}$ is in $\mathcal{V}_f$, then for all agents $i\in\mathcal{A}$ with $\{i,j\}\in\bar{\mathcal{E}}$, we have $i\in\mathcal{V}_f$. By contrapositive, if an agent $i\in\mathcal{A}$ with $\{i,j\}\in\bar{\mathcal{E}}$ is not in $\mathcal{V}_f$, then task $j\in\mathcal{B}$ is not in $\mathcal{V}_f$. Similarly, if an agent $i\in\mathcal{A}$ is in $\mathcal{V}_f$, then $m_i\in\mathcal{V}_f$, and if $m_i\notin \mathcal{V}_f$, then agent $i\in\mathcal{A}$ is not in $\mathcal{V}_f$.

Assume for contradiction that there exists $i\in\mathcal{A}$ such that there exists an alternating path $\mathcal{P}$ between $i$ and $m_{a_n}$ and that $i\notin\mathcal{V}_f$. Let $\mathcal{P}:=\{\{v_k,v_k+1\}\in\bar{\mathcal{E}}|k=0,1,...,K, v_0=i, v_{K+1}=m_{a_n}\}$. We have that $v_0=i\notin \mathcal{V}_f$. From the arguments above, if $v_k\notin \mathcal{V}_f$, then $v_{k+1}\notin\mathcal{V}_f$. By applying this to all edges in $\mathcal{P}$, we have that $v_{K+1}=m_{a_n}\notin \mathcal{V}_f$. This is a contradiction, as $m_{a_n}$ is the root of the tree. Therefore, $i$ must be in $\mathcal{V}_t$.

\subsection{Proof of Lemma~\ref{lem:augpath}} \label{ap:lem:augpath}

At every iteration of the while-loop in Function 1, there exists an alternating path $\mathcal{P}$ between the current vertex $t$ and root $\bar{j}$. Consider the set of agents incident to edges in $\mathcal{P}$, i.e., $K_\mathcal{P}:=\{i\in\mathcal{A}|\{i,j\}\in \mathcal{P}\}$, and consider function $g:\mathcal{A}\rightarrow \mathcal{B}$ mapping an agent in the tree to its parent vertex. From Lines~\ref{line:nudone} and~\ref{line:nufor}, for all agents $i\in K_\mathcal{P}$, $\nu_i=g(i)$. From Lines~\ref{line:nuinit} and~\ref{line:nuback}, we have that $\nu_i=m_i$, for all $i\notin K_\mathcal{P}$. Thus, at every iteration of the while-loop we have $\{\{i,\nu_i\}|i\in\mathcal{A}\}=\{\{i,m_i\}|i\notin K_\mathcal{P}\}\cup \{\{i,g(i)\}|i\in K_P\}$.

From Line~\ref{line:finish}, we get $\{\{i,\nu_i\}|i\in\mathcal{A},\nu_i\neq -1\}=\{\{i,m_i\}|i\notin K_\mathcal{P},m_i\neq-1\}\cup \{\{i,g(i)\}|i\in K_\mathcal{P}\}=\bar{\mathcal{M}}\oplus \mathcal{P}$. Thus, if $\mathcal{P}$ is an augmenting path, then the function returns $\mathcal{M}_{\nu}=\bar{\mathcal{M}}\oplus \mathcal{P}$ as required. On the other hand, if no augmenting path exists, the search terminates at $t=\bar{j}$, i.e., $\mathcal{P}=\emptyset$, $K_\mathcal{P}=\emptyset$, and $\{\{i,\nu_i\}|i\in\mathcal{A},\nu_i\neq -1\}$$=$$\{\{i,m_i\}|i\in\mathcal{A},m_i\neq -1\}=\bar{\mathcal{M}}$.

\subsection{Proof of Theorem~\ref{thm:pruneBAP}} \label{ap:thm:pruneBAP}

From Proposition~\ref{lem:algcorrect}, pruneBAP returns a bottleneck assignment of $\mathcal{G}_b$. From Lemmas~\ref{lem:maxconsensus} and~\ref{lem:augpath}, the functions $\texttt{\textsc{MaxEdge}}()$ and $\texttt{\textsc{AugPath}}()$ in pruneBAP can be implemented such that they satisfy Assumptions~\ref{as:distedge} and~\ref{as:distcom}. Edge removal satisfying Assumptions~\ref{as:distedge} and~\ref{as:distcom} can be implemented as shown in~(\ref{eq:distphi}),~(\ref{eq:distremedge}), and~(\ref{eq:distremmatch}).

\subsection{Proof of Proposition~\ref{prop:worst-case}} \label{ap:prop:worst-case}

In each iteration of the while-loop of pruneBAP, at least one edge in $\mathcal{E}_b$ is removed. There are therefore at most $|\mathcal{E}_b|=mn$ iterations of the while-loop. The while-loop of pruneBAP itself contains two components that depend on time steps. The function $\texttt{\textsc{MaxEdge}}()$ requires at most $D$ time steps for completion. The distributed search for an augmenting path has order $\mathcal{O}(nD)$, where each iteration of the while-loop in Function~\ref{alg:augpath} requires $D$ time steps and there are at most $2n-1$ iterations.

\subsection{Proof of Lemma~\ref{lem:BFS}} \label{ap:lem:BFS}

Function~\ref{alg:bfs} extends Function~\ref{alg:augpath} by exploring multiple branches of the alternating tree in parallel. In Function~\ref{alg:augpath}, only one task $t$ is denoted as the current vertex per iteration of the while loop; in Function~\ref{alg:bfs}, a set of tasks $\bar{\mathcal{B}}$ as considered as the current vertices.	At every iteration of the while loop in Function~\ref{alg:bfs}, for all $t\in \bar{\mathcal{B}}$, there exists an alternating path $P_t$ between $t$ and $\bar{j}$. This holds by construction, since an agent $i\in F$ has its matched task $m_i$ inserted into the set $\bar{\mathcal{B}}$ only after $i$ is explored. We alternate between edges in the matching and not in the matching as outlined in Definition~\ref{as:tree}. If one or more unmatched agents are explored in Line~\ref{line:barS}, then the search is successful and it remains for one of the augmenting paths to be selected. Following Proposition~\ref{prop:systematic}, no augmenting path exists if $\bar{S}$ is empty. In such a case, all agents $i\in\mathcal{A}$ for which there exists an alternating path between $i$ and root $\bar{j}$ have been explored and none of these agents are free vertices.

\subsection{Proof of Lemma~\ref{thm:bound}} \label{ap:thm:bound}

By definition, $w(e_3)=\max_{e\in\mathcal{M}_3}w(e)\leq \max_{e\in\mathcal{M}}w(e)$ for any arbitrary $\mathcal{M}\in\mathcal{C}(\mathcal{G}_3)$. Since $\mathcal{M}'=\mathcal{M}_1\cup\mathcal{M}_2\in\mathcal{C}(\mathcal{G}_3)$, we have $w(e_3)\leq \max_{e\in\mathcal{M}'}w(e)=\max\{w(e_1),w(e_2)\}$.

\subsection{Proof of Proposition~\ref{prop:cluster}} \label{ap:prop:cluster}

Path $\mathcal{P}=\{e_c\}$ is trivially an alternating path relative to $\mathcal{M}$ and edge $e_c\in \mathcal{M}$, therefore $\mathcal{P}\subseteq \phi(\mathcal{G}_b,\mathcal{M})$. It remains to show that there does not exist another alternating path $\mathcal{P}'\subseteq \phi(\mathcal{G}_b,\mathcal{M})$ relative to $\mathcal{M}$ between $a_c$ and $b_c$. To this end, assume for contradiction that aside from the path $\mathcal{P}=\{e_c\}$ there exists another alternating path $\mathcal{P}'$ relative to $\mathcal{M}$ between $a_c$ and $b_c$, where $\mathcal{P}'\neq \mathcal{P}$ and $\mathcal{P}'\subseteq \phi(\mathcal{G}_b,\mathcal{M})$. From Definition~\ref{def:path}, a path is constructed from a sequence of distinct vertices, and $a_c$ and $b_c$ are the first and last vertices of the sequence so $\mathcal{P}'$ cannot contain $e_c$ without repeating $a_c$ or $b_c$ in the sequence. Since $\mathcal{P}'\neq \mathcal{P}$, it is guaranteed that $e_c\notin \mathcal{P}'$ and $\mathcal{P}' \subseteq \phi(\mathcal{G}_b,\mathcal{M})\backslash\{e_c\}$. It follows that $\mathcal{P}'$ is an augmenting path relative to $\mathcal{M}\backslash \{e_c\}$. This contradicts the assumption that $e_c$ is a critical bottleneck edge of $\mathcal{G}_b$ relative to $\mathcal{M}$.

\subsection{Proof of Lemma~\ref{lem:reduce}} \label{ap:lem:reduce}

Without loss of generality, let $e_1=\{a_1,b_1\}$. By Proposition~\ref{prop:cluster}, $e_1$ is the only alternating path between $a_1$ and $b_1$ in edge set $\phi(\mathcal{G}_1,\mathcal{M}_1)$. Assume there does not exist vertices $i$ and $j$ such that all i., ii., and iii. are true. Thus, $e_1$ is the only alternating path between $a_1$ and $b_1$ in $\phi(\mathcal{G}_3,\mathcal{M}_1\cup\mathcal{M}_2)$. By Definition~\ref{def:crit}, $e_1$ is also a critical bottleneck edge of $\mathcal{G}_3$ since $e_1\in\arg\max_{e\in\mathcal{M}_1\cup\mathcal{M}_2} w(e)$ and there does not exist an augmenting path in $\phi(\mathcal{G}_3,\mathcal{M}_1\cup\mathcal{M}_2)\backslash\{e_1\}$ relative to $(\mathcal{M}_1\cup\mathcal{M}_2)\backslash\{e_1\}$. Thus, $w(e_3)=w(e_1)$.

\subsection{Proof of Theorem~\ref{thm:reduce2}} \label{ap:thm:reduce2}

The necessary condition for $w(e_3)<w(e_1)$ holds from Lemma~\ref{lem:reduce}. We now prove the sufficient condition. Assume conditions i., ii., and iii. hold. Then, aside from path $\mathcal{P}=\{e_1=\{a_1,b_1\}\}$, there exists another alternating path $\mathcal{P}'$ between $a_1$ and $b_1$, which does not contain the edge $e_1$. Namely, the alternating path $\mathcal{P}'$ constructed from the union of the alternating paths between $a_1$ and $b'$, $b'$ and $i$, $i$ and $j$, $j$ and $a'$, and $a'$ and $b_1$. Thus, there exists an augmenting path $\mathcal{P}'\subseteq\phi(\mathcal{G}_3,\mathcal{M}_1\cup\mathcal{M}_2)\backslash\{e_1\}$ relative to $(\mathcal{M}_1\cup\mathcal{M}_2)\backslash\{e_1\}$. From Proposition~\ref{lem:distBAP}, there exists an MCM $\mathcal{M}'$ of $\mathcal{G}_3$ such that $\mathcal{M}'\subseteq \phi(\mathcal{G}_b,\mathcal{M}_1\cup\mathcal{M}_2)\backslash\{e_1\}$. By the assumptions on $e_1$, $\phi(\mathcal{G}_3,\mathcal{M}_1\cup\mathcal{M}_2)\backslash\{e_1\}$ contains only edges with weights strictly smaller than $w(e_1)$. Thus, there exists an MCM of $\mathcal{G}_3$ with all edges having weight smaller than $w(e_1)$, therefore, $w(e_3)$ must be smaller than $w(e_1)$.

\end{document}